\documentclass[11pt, letterpaper,fleqn]{amsart}
\usepackage{amssymb,amsmath,breqn}
\usepackage[foot]{amsaddr}
\usepackage{mathtools} 
\usepackage{hyperref}
\usepackage{blkarray} 
\usepackage{graphicx}

\newcommand{\RNum}[1]{\uppercase\expandafter{\romannumeral #1\relax}}

\usepackage{comment}
\usepackage{cleveref}

\usepackage{enumerate}

\newcommand{\ba}{\bar{a}}
\newcommand{\bb}{\bar{b}}
\newcommand{\bi}{\bar{i}}
\newcommand{\bj}{\bar{j}}
\newcommand{\bk}{\bar{k}}
\newcommand{\bl}{\bar{l}}
\newcommand{\bm}{\bar{m}}
\newcommand{\bn}{\bar{n}}
\newcommand{\bp}{\bar{p}}
\newcommand{\bq}{\bar{q}}

\newcommand{\bz}{\bar{z}}

\newcommand{\bI}{\bar{I}}
\newcommand{\bJ}{\bar{J}}
\newcommand{\bK}{\bar{K}}
\newcommand{\bL}{\bar{L}}

\newcommand{\bM}{\bar{M}}

\newcommand{\bpartial}{\bar{\partial}}
\newcommand{\pbp}{\partial \bar{\partial}}

\newcommand{\fI}{\mathfrak{I}}

\newcommand{\fRe}{\mathfrak{Re}}

\newcommand{\cC}{\mathcal{C}}

\newcommand{\cH}{\mathcal{H}}

\newcommand{\cL}{\mathcal{L}}

\newcommand{\cP}{\mathcal{P}}
\newcommand{\cT}{\mathcal{T}}

\newcommand{\tr}{\mbox{tr}}

\newcommand{\fJ}{\mathfrak{J}}

\usepackage{xcolor}

\newcommand{\ol}{\overline}
\newcommand{\ul}{\underline}

\newcommand{\p}{\partial}
\newcommand{\vp}{\varphi}

\newtheorem{theorem}{Theorem}[section]
\newtheorem{lemma}[theorem]{Lemma}

\newtheorem{proposition}[theorem]{Proposition}

 \theoremstyle{definition}
\newtheorem{definition}[theorem]{Definition}

\abovedisplayskip
\belowdisplayskip

\theoremstyle{remark}

\numberwithin{equation}{section}

\begin{document}

\title[Complex Monge-Amp\`ere equation for positive $(p,p)$-forms]
{Complex Monge-Amp\`ere equation for positive $(p,p)$-forms on compact K\"ahler manifolds}


\author{Mathew George}
\address{Department of Mathematics, Purdue University,
         West Lafayette, IN 47907, USA}
\email{georg233@purdue.edu}

\date{}

\begin{abstract}

A complex Monge-Amp\`ere equation for differential $(p,p)$-forms is introduced on compact K\"ahler manifolds. For any $1 \leq p < n$, we show the existence of smooth solutions unique up to adding constants. For $p=1$, this corresponds to the Calabi-Yau theorem proved by S. T. Yau, and for $p=n-1$, this gives the Monge-Amp\`ere equation for $(n-1)$ plurisubharmonic functions studied by Tosatti-Weinkove. For other $p$ values, this defines a non-linear PDE that falls outside of the general framework of Caffarelli-Nirenberg-Spruck. Further, we define a geometric flow for higher-order forms that extends the K\"ahler-Ricci flow naturally to $(p,p)$-forms. As a consequence of our main theorem, we show that this flow exists in a maximal time interval and can be shown to converge under some assumptions. A modified flow is introduced and the convergence of the associated normalized flow is shown.



{\em Mathematical Subject Classification (MSC2020):}
35J15, 35J60, 58J05.
\end{abstract}

\maketitle

\section{Introduction}

\medskip

The study of fully nonlinear PDEs on K\"ahler manifolds has led to many important developments over the last 50 years. Many of these equations are related to the complex Monge-Amp\`ere equation, which plays a fundamental role in this theory. In this paper, we introduce a new class of equations for $(p,p)$-forms for $1 \leq p \leq n-1$ that closely resembles the complex Monge-Amp\`ere equation. 

 Given a compact K\"ahler manifold $(M, \omega)$, let $\Omega$ be a real, positive differential $(p,p)$-form on $M$. Then

$$\Omega_{u } = \Omega + \sqrt{-1} \pbp u \wedge \omega^{p-1} $$

\noindent defines another real $(p,p)$-form. If $\Omega$ is closed, then $\Omega_u$ is a closed form such that $[\Omega_u] = [\Omega]$ in $H^{p,p}(M, \mathbb R)$. For a smooth function $F$ on $M$, we define the equation

\begin{equation}\label{pde1}
    \mu_{\Omega_u} = e^{F} \omega^n
\end{equation}

\noindent where $\mu_{\Omega_u}$ is the volume form associated to $\Omega_u$. In local coordinates, $\mu_{\Omega}$ is given by $(\det{(\Omega_{I \bJ})})^{\frac{n}{pN}}$, where $\Omega_{I \bJ}$ is the matrix obtained from the coefficients of the $(p,p)$-form $\Omega$ using lexicographic ordering on its indices, and $N = {n \choose p}$. See Section $2$ for precise definitions. We prove the following.

\begin{theorem}\label{main-theorem}
    Given a compact K\"ahler manifold $(M, \omega)$, a smooth function $F$ on $M$, and a positive $(p,p)$-form $\Omega$, there exists a unique smooth function $u$ with $\sup\limits_M u = 0$, and a unique constant $b$ such that

    \begin{equation}\label{main-pde}
        \mu_{\Omega_u} = e^{F + b} \omega^n
    \end{equation}

    \

    \noindent where $\Omega_u = \Omega + \sqrt{-1} \pbp u \wedge \omega^{p-1} $ is a positive form.
\end{theorem}

In 1978, S.T. Yau \cite{Yau78} solved the complex Monge-Amp\`ere equation and used it to show the existence of K\"ahler metrics with prescribed volume forms on compact K\"ahler manifolds. Since for K\"ahler manifolds, the Ricci curvature is given by second-order variations of $\log$ of the volume form, this proves the existence of Ricci-flat K\"ahler metrics on compact K\"ahler manifolds with $c_1(M) = 0$. More precisely, Yau's theorem shows that given a $(1,1)$-form $\Psi$ in $c_1(M)$ and a K\"ahler class $[\omega]$, there exists a unique K\"ahler metric $\omega'$ in $[\omega]$ such that 

$$\operatorname{Ric}(\omega')= \Psi.$$

\

Similar to this, one could define a ``Ricci form" for positive $(p,p)$-forms by considering $\operatorname{Ric}(\Omega) = -\sqrt{-1} \pbp \log{(\mu_{\Omega})}$. Then equation \eqref{main-pde} can equivalently be written as $\operatorname{Ric}(\Omega_u) = \Psi$, with $\Psi \in [\operatorname{Ric}(\tilde{\Omega})] $ for some positive form $\tilde{\Omega}$. This could have a natural interpretation, since $\mu_{\Omega}$ coincides with the volume form associated to the metric $\omega$, if $\Omega = \omega^{p}$. It can also be shown that the cohomology class $[\operatorname{Ric}(\Omega)]$ is independent of $\Omega$.


By the works of Cherrier \cite{Cherrier87}, Guan-Li \cite{GL10}, Tosatti-Weinkove \cite{TW10}, the estimates for complex Monge-Amp\`ere equation was extended to compact Hermitian manifolds. Fu-Wang-Wu \cite{FWW10, FWW15} first studied the $(n-1)$ form-type Monge-Amp\`ere equation, and later in \cite{TW17}, Tosatti-Weinkove solved the Monge-Amp\`ere equation for $(n-1)$ plurisubharmonic functions on K\"ahler manifolds. This was also extended to Hermitian manifolds by Tosatti-Weinkove in \cite{TW19}. The positive $(n-1,n-1)$ form

$$\omega_u^{n-1} = \omega^{n-1} + \sqrt{-1} \pbp u \wedge \omega^{n-2}$$

\noindent defines a positive $(1,1)$-form ${\omega}_u$ uniquely. Hence one can write the $(n-1)$ plurisubharmonic Monge-Amp\`ere equation as

$${\omega}_u^n = e^F \omega^n,$$

\noindent since $\omega_u$ is known. This is not the case here as the Hodge-star operator applied to a $(p,p)$-form gives an $(n-p, n-p)$ form. So instead we work directly with the $(p,p)$-form using multi-indices.

An important detail is that equation \eqref{pde1} is not symmetric in the eigenvalues of the Hessian of $u$. Although it is concave in $D^2u$, it need not satisfy some familiar properties \cite{CNS85} such as $ F^{i \bi} \leq F^{j \bj}$ for $\lambda_i \geq \lambda_j$, that is true for Monge-Amp\`ere type equations $F(\omega + \sqrt{-1} \pbp u) = \psi$.

The main challenge in the proof of Theorem \ref{main-theorem} is obtaining a $C^2$ a priori estimate for the solution $u$. We briefly outline the difficulties and some of the ideas involved in the proof. In the maximum principle argument for $C^2$ estimates, the third-order terms that are obtained by differentiating the PDE twice must be used to control the non-positive third-order terms coming from the test function. Given the more complicated structure of this equation, there are some difficulties in obtaining a direct comparison between them. We overcome this by splitting the multi-indices into cases, when $Z_{I \bI}$ (the matrix of $\Omega_u$) is close to the trace of $(Z_{I \bJ})$ and when it is much smaller. The part when it is smaller is controlled by using the concavity terms from the second derivative of the equation. The other terms can be handled using several arguments involving the smallness of the inverse terms $Z^{I \bI}$, and some larger positive terms involving $\lambda_1^2$ that originate from choosing the test function appropriately, similar to Hou-Ma-Wu \cite{HMW10}. The second difficulty is that there is no way of diagonalizing the matrix of a general $(p,p)$-form by choosing a coordinate system. This is dealt with by using the components of the eigenvectors of $(Z_{I \bJ})$, and obtaining the right estimates for them with respect to the trace.





In the recent years, a vast array of new techniques have been developed for fully non-linear equations on complex manifolds, most notably by the works of Tossatti-Weinkove \cite{TW10,TW17,TW19}, Szekelyhidi \cite{Szekelyhidi18}, Chu-Tosatti-Weinkove \cite{CTW17}, Dinew-Kolodziej \cite{DK17}, Phong-Picard-Zhang \cite{PPZ21}, Guo-Phong-Tong \cite{GPT23, GP23}, Hou-Ma-Wu \cite{HMW10}, George-Guan-Qiu \cite{GGQ22}, Guan \cite{Guan14} among others. But there are still many difficulties in solving nonlinear PDEs involving differential forms that are important in many areas of geometry and mathematical physics. It is possible that the methods introduced here could find wider applicability. An extension to fully nonlinear case is considered in \cite{GG24}, where other symmetric functions of eigenvalues of $\Omega_u$, satisfying a rank condition, are considered.

\

 In the latter part of the paper, we turn our attention to the associated parabolic equation to study an evolution of positive forms. For $F \in C^{\infty}(M)$ and $\tilde{\Omega}$ a positive form, consider the scalar parabolic equation

\begin{equation}\label{scalar-flow}
    \frac{\p u}{ \p t} = \frac{pN}{n}\left(\log{\frac{\mu_{\Omega_u}}{\mu_{\tilde{\Omega}}}}  - F\right)
\end{equation}

\noindent with $u(0) = 0$. Here $\Omega_u = \Omega_0 + \sqrt{-1} \pbp u \wedge \omega^{p-1}$. Then we prove the following theorem.

\begin{theorem}\label{scalar-flow-theorem}
Let $(M, \omega)$ be a compact K\"ahler manifold, and $\Omega_0$ a positive form. Then there is a unique smooth solution $u$ to the equation \eqref{scalar-flow} for $t \in [0, \infty)$, satisfying $\Omega_u >0$ for $t \in [0, \infty)$. In addition, the normalized solution
$$\tilde{u} = u -\frac{1}{V}\int_M u \omega^n$$
converges smoothly to a solution $u_{\infty}$ of the elliptic equation $\mu_{\Omega_{u_{\infty}}} = e^{F+b} \mu_{\tilde{\Omega}} $ for some constant $b = b_{\infty}$.
\end{theorem}

Some details of the proof of this theorem are omitted in this paper, since it is similar to the elliptic estimate, and the convergence is given in detail in other works. We instead outline the argument for a normalized version of this flow in Theorem \ref{normalized-estimate-theorem}. The parabolic case of a general fully-nonlinear equation of real forms is dealt with in \cite{GG24}. Also see \cite{Cao85}, \cite{Gill11}, \cite{PT21}, and \cite{George21}.

Consider a flow of $(p,p)$-forms given by

\begin{equation}\label{flow1}
\begin{aligned}
    &\frac{\p}{ \p t}  \Omega(z, t) =   \frac{pN}{n}\sqrt{-1} \pbp \log{\mu_{\Omega(z, t)}} \wedge \omega^{p-1} ,\\
     &\Omega(z, 0) = \Omega_0(z) \in \Lambda^{p,p}M
    \end{aligned}
\end{equation}

\noindent for a fixed $(p,p)$-form ${\Omega}_0$. This reduces to the K\"ahler-Ricci flow for $p=1$, and for $p= n-1$ such an equation was studied by Gill \cite{Gill14}, and Tosatti-Weinkove \cite{TW17}. 

\

Denote $\operatorname{Ric_p}(\Omega) = -\sqrt{-1} \pbp \log{\mu_{\Omega}} \wedge \omega^{p-1}$. As a corollary to the Theorem \ref{scalar-flow-theorem},

\begin{theorem}\label{existence-time-theorem}
    Given any $\Omega_0 \in \Lambda^{p,p}(M, \mathbb R)$ on a compact K\"ahler manifold $(M, \omega)$, there exists a smooth positive solution $\Omega_u$ to \eqref{flow1}  for $t \in [0, T)$ given by 
\begin{equation}\label{form-flow-solution}
    \Omega_u = \Omega_0  - \frac{pN}{n}t \operatorname{Ric_p}(\tilde{\Omega})+ \sqrt{-1} \pbp u \wedge \omega^{p-1},
\end{equation}

\noindent where 
   \begin{equation}
       \begin{aligned}
             T = \sup\{&t \in (0, \infty): \; \text{$\exists$ a smooth function $\tilde{u}$, and a positive form $\tilde{\Omega}$ such that }\\& \Omega_0 - t \frac{pN}{n}\operatorname{Ric_p}(\tilde{\Omega}) + \sqrt{-1} \pbp \tilde{u} \wedge \omega^{p-1} >0\}.
       \end{aligned}
   \end{equation}

  
  \end{theorem}

\

If we also make the assumption that $\Omega_0$ is positive, then the following holds.

\begin{theorem}
    Assume that $\Omega_0$, and $\tilde{\Omega}$ are real positive $(p,p)$-forms. Then the equation
    \begin{equation}\label{flow2}
\begin{aligned}
    &\frac{\p}{ \p t}  \Omega =   \frac{pN}{n}\sqrt{-1} \pbp \log\frac{{\mu_{\Omega}}}{\mu_{\tilde{\Omega}}} \wedge \omega^{p-1} ,\\
     &\Omega(z, 0) = \Omega_0(z) \in \Lambda^{p,p}M,
    \end{aligned}
\end{equation}

\noindent has a solution in $[0, \infty)$ given by 

\begin{equation}\label{form-solution}
    \Omega_u = \Omega_0 + \sqrt{-1} \pbp u \wedge \omega^{p-1},
\end{equation}

\noindent where $u$ satisfies \eqref{scalar-flow} with $F= \log(\mu_{\tilde{\Omega}/\omega^n})$. In addition if $[\Omega_0] \in H^{p,p}(M, \mathbb R)$, then $[\Omega_u]\in H^{p,p}(M, \mathbb R)$, and the flow converges as $t \to \infty$, so that

  \begin{equation}\label{form-infinite}
       \Omega_{\infty} = \Omega_0 + \sqrt{-1} \pbp u_{\infty} \wedge \omega^{p-1}.
  \end{equation}
 
\end{theorem}

From \eqref{form-solution} it is clear that the flow preserves the cohomology class. It is natural to suspect that it converges to a canonical representative in this class.





Motivate by this we consider 

\begin{equation}\label{flow4}
\begin{aligned}
   &\frac{\partial \Omega}{\partial t} =\frac{pN}{n} \sqrt{-1} \pbp \log{\frac{ \mu_{\Omega}}{\mu_{\tilde{\Omega}}}}\wedge \omega^{p-1} + \Omega_0,\\
   &\Omega|_{t=0} = \Omega_0 \in H^{p,p}(M, \mathbb R)
\end{aligned}
\end{equation}

A solution could be written of the form

\begin{equation}
    \Omega_u = (1 + t)\; \Omega_0 + \sqrt{-1} \pbp u \wedge \omega^{p-1}
\end{equation}

\noindent for a $u$ satisfying 

\begin{equation}\label{parabolic-ma-eqn1}
    \frac{\partial u}{\partial t} = 
    \frac{pN}{n}\log\frac{\mu_{\Omega_u}}{\mu_{\tilde{\Omega}}} 
\end{equation}

\noindent with $u(0) = 0$. The form $\Omega_0$ being positive is enough to show the long-time existence of the flow. By a reparametrization of \eqref{flow4}, we get

\begin{equation}\label{normalized-form-flow1}
\begin{aligned}
 &\frac{\partial\Omega}{\partial t} = \frac{pN}{n} \sqrt{-1} \pbp \log{\frac{ \mu_{\Omega}} {\mu_{\tilde{\Omega}}}}\wedge \omega^{p-1}+ \Omega_0 - \Omega\\
   &\Omega|_{t=0} = \Omega_0 \in H^{p,p}(M, \mathbb R)
\end{aligned}
\end{equation}

A solution to \eqref{normalized-form-flow1} is given by 

\begin{equation}
    \Omega_u = \Omega_0 + \sqrt{-1} \pbp u \wedge \omega^{p-1}
\end{equation}

\noindent where $u$ satisfies the normalized $(p,p)$-form parabolic Monge-Amp\`ere equation

\begin{equation}\label{normalized-ma-eqn1}
\begin{aligned}
    & \frac{\partial u}{\partial t}  = \frac{pN}{n}\log\frac{\mu_{\Omega_u}}{\mu_{\tilde{\Omega}}} - u, \text{ on } M\\
   & u(0) = 0.
\end{aligned}
\end{equation}




\begin{theorem}
    Let $(M, \omega)$ be a compact K\"ahler manifold with $\Omega_0 \in H^{p,p}(M, \mathbb R)$ a positive form. Then there exists a positive $(p,p)$-form $\Omega_u$ that satisfies \eqref{normalized-form-flow1} for $t$ in $[0, \infty)$, given by 

    \begin{equation}
    \Omega_u = \Omega_0 + \sqrt{-1} \pbp u \wedge \omega^{p-1}.
\end{equation}

\

\noindent for a solution $u$ of \eqref{normalized-ma-eqn1}. In addition, as $t \to \infty$, $\Omega_u$ converges smoothly to $\Omega_{\infty}$ which satisfies

\begin{equation}
   \Omega_{\infty} = \Omega_0 + \frac{pN}{n}\sqrt{-1} \pbp \log{\frac{ \mu_{\Omega_{\infty}}}{\mu_{\tilde{\Omega}}}}\wedge \omega^{p-1}
\end{equation}

\end{theorem}




\

More details of the flow are given in Section \ref{section6}.

 \bigskip

\section{Preliminaries}

\subsection{\bf Matrix associated to a $(p,p)$-form:}
For a fixed number $p \leq n$, define the set of multi-indices of size $p$,
$$\fI_p = \{(i_1, \hdots, i_p  ) : 1 \leq i_1 < \hdots < i_p \leq n,\; i_j \in \mathbb N\}$$

Clearly there are $N= \dfrac{n!}{p! (n-p)!}$ elements in $\fI_p$.
We now fix the order in $\fI_p$ {\sl lexicographically}.
\[ I = (i_1, \ldots, i_p) < J = (j_1, \ldots, j_p) \]
if $i_l < j_l$ for the first non-equal pair. So we write $ \fI_p = \{I_1, \ldots, I_N: I_i < I_j \; \mbox{if $i < j$}\} \;$.

Given a $(p,p)$-form $\Omega$ in local coordinates, define a matrix associated to it as follows. For $I  = (i_1, \ldots, i_p) \in \fI_p$ write 
\[ dz_I = dz_{i_1} \wedge \cdots \wedge dz_{i_p}, \;\;
 d\bz_I = d \bz_{i_1} \wedge \cdots \wedge d \bz_{i_p}. \]

A $(p,p)$-form in local coordinates can be written as
$$\Omega = p! (\sqrt{-1})^{p^2} \Omega_{I \bJ} dz_I \wedge d \bz_{J}$$

Then following the lexicographic ordering of multi-indices, the components $\Omega_{I \bJ}$ are arranged as the entries of a matrix. We will denote the matrix associated to $\Omega$ by the same notation
$\Omega = (\Omega_{I \bJ})$. For the $(1,1)$-form $\omega = \sqrt{-1} g_{i \bj} dz_i \wedge d \bz_j$, 

$$\omega^p = p! (\sqrt{-1})^{p^2} g_{I \bJ} dz_{I} \wedge d\bz_{J}$$

\noindent where $g_{I \bJ} = \det{(g_{i \bj})_{i \in I, \; j \in J}}$.

\subsection{\bf Notations and Conventions:}

We will use the following notation that is useful in the study of this equation. For $I' = (i_1, \hdots, i_{p-1})$ in $\fI_{p-1}$, denote for any $k \notin I'$, and $i_{j} < k < i_{j+1}$,

$$ I'_k = (i_1, \hdots , i_{j}, k , i_{j+1}, \hdots, i_{p-1}),$$

\noindent so that $I'_k \in \fI_p$. The notation $(-1)^{(i| I)}$ would stand for the sign of the permutation so that

$$dz_{i} \wedge d z_{I'} = (-1)^{(i| I)} dz_{I'_i}.$$

In local coordinates, we write

$$\Omega_{u} = (p-1)! (\sqrt{-1})^{p^2}Z_{I \bJ} dz_I \wedge dz_{\bJ},$$

\noindent and $(Z^{I \bJ})$ is the inverse matrix of $(Z_{I \bJ})$. For simplicity we will assume the following summation convention unless otherwise stated.

$$(-1)^{(i|I'_i)+(j|I'_j)}Z^{I'_i \bI'_j} \xi_i \xi_{\bj}$$

\noindent would stand for the summation

$$\sum\limits_{I' \in \fI_{p-1}} \sum\limits_{i,j \notin I'} (-1)^{(i|I'_i)+(j|I'_j)}Z^{I'_i \bI'_j} \xi_i \xi_{\bj}.$$

If there is a deviation from this notation, for example if $I'_i$ is in a particular subset, that will be made clear.
Define the positive cone in $\mathbb R^N$ by

$$\Gamma_N = \{\Lambda = (\Lambda_1, \hdots, \Lambda_N) \in \mathbb R^N: \Lambda_i > 0\; \text{ for } 1\leq i \leq N \}.$$

For any subset $S$, $S^c$ will denote its complement.

\subsection{\bf Volume form of a positive $(p,p)$-form:} A $(p,p)$-form $\Omega = p! (\sqrt{-1})^{p^2} \Omega_{I \bJ} dz_I \wedge d \bz_{J}$ is called positive if the matrix associated to it is positive-definite. In addition, one can define the eigenvalues of $\Omega$ as roots of the polynomial $$\det{(\Omega - t \omega^p)} = 0.
$$

We could also define this as the eigenvalues of the linear transformation on $\Lambda^p \; TM$, obtained by raising the anti-holomorphic indices of $\Omega$ by the metric $g$.

We will denote the eigenvalues of $\Omega$ by $\Lambda_N(\Omega) \leq \hdots \leq \Lambda_1(\Omega).$

\begin{definition}
Let $\Omega$ be a positive $(p,p)$ from. The associated {\sl volume form} of 
$\Omega$ is defined locally by 
\[ \mu_{\Omega} = (\sqrt{-1})^{n} \det (\Omega_{I\bJ})^{\frac{n}{p N}} dz_1 \wedge d\bz_1 \wedge \cdots \wedge dz_n \wedge d\bz_n. \]
\end{definition}

 In particular, if $\Omega = \eta^{p}$ for a positive $(1,1)$-form $\eta$, then by Sylvester-Franke theorem
 \[ \mu_{\Omega} = \frac{\eta^n}{n!} . \]

It is not difficult to check that these definitions are invariant under coordinate changes. Indeed, if $A$ is the change of coordinates matrix for $(1,1)$-forms at a point, then the matrix corresponding to change of coordinates of $(\Omega_{I \bJ})$ is given by $C_p(A) = \wedge^p A$, the $p^{th}$ compound matrix associated to $A$. It is also the case that $C_p(A)$ is functorial in $A$, and 

$$\det(C_p(A))^{\frac{n}{pN}} = \det(A).$$

Hence $\mu_{\Omega}$ defines a global volume form on $M$. We now prove that positive $(p,p)$-forms as defined above are also positive as currents given by differential forms. 

\begin{proposition}\label{positive-theorem}
    Let $M$ be a compact complex manifold. Then a positive $(p,p)$-form is also positive as a current.
\end{proposition}

\begin{proof} Using the convention in \cite{Demailly12}, a positive current $T$ is defined by $\langle T, \beta\rangle >0$, for any strongly positive form $\beta$.

    Let $Z = \sqrt{-1}^{p^2} Z_{I \bJ} dz_I \wedge d \bz_{J}$ be a positive $(p,p)$-form as defined above. Then a strongly positive $(n-p,n-p)$ form is a linear combination with positive coefficients of forms of the type $\beta =\sqrt{-1}^{(n-p)^2} \alpha \wedge \bar{\alpha}$, where
    $$\alpha = \alpha_1 \wedge \hdots \wedge \alpha_{n-p} ,$$ 

\noindent $\alpha_i = \sum_j\alpha_i^{j} dz_j$ are elementary $1$-forms. For a multi-index $K= (k_1,\hdots, k_{n-p})$ let 

$$\alpha_{K} =\sum_{\tau \in \sigma_p}sgn(\tau)\alpha^{k_1}_{\tau(1)}\hdots \alpha_{\tau(n-p)}^{k_{(n-p)}}$$
\ 

 Here $\sigma_p$ is the set of all permutations of $\{1, \hdots ,n-p\}$. We write $\alpha = \sum_K \alpha_K dz_K.$ Define the sign $(-1)^{(K|K^c)}$ by

$$dz_{K }\wedge dz_{K^c} = (-1)^{(K|K^c)} dz_1 \wedge dz_2 \wedge \hdots \wedge dz_n. $$

Then

\begin{equation}
\begin{aligned}
    \frac{Z \wedge \beta}{d\mathrm{vol}} &=
    \sum_{I, J \in \fI_p} Z_{I \bJ} (-1)^{(I|I^c)}\alpha_{I^c}\ol{(-1)^{(J|J^c)}\alpha_{J^c}}
 > 0 
    \end{aligned}
\end{equation}

\noindent positivity of which can be seen by taking a complex vector $\xi_I = (-1)^{(I|I^c)}\alpha_{I^c}$. 
This implies that the pairing 
$$\langle Z, \beta \rangle = \int_M Z \wedge \beta > 0$$

\noindent for any non-zero strongly positive form $\beta$.




\end{proof}

\subsection{\bf Curvature and connection:} In local coordinates $z = (z_{1}, \ldots,z_{n})$, we shall write 
\begin{equation}
\label{gqy-pre213}
\omega = \sqrt{-1} g_{i\bj} dz_i \wedge dz_{\bj}
\end{equation}
and $\{g^{i\bj}\} = \{g_{i\bj}\}^{-1}$.
The Christoffel symbols $\Gamma_{ij}^k$ is defined by
\[ \nabla_{\frac{\partial}{\partial z_i}} \frac{\partial }{\partial z_j}
= \Gamma_{ij}^k \frac{\partial}{\partial z_k}, 
\;\;  \]
where $\nabla$ is the Levi-Civita connection on $(M, \omega)$, so 
\[ \Gamma_{ij}^k = g^{k\bl} \frac{\partial g_{j\bl}}{\partial z_i}. \]

For commuting covariant derivatives of a function, we will use the following commutation formulas.

\begin{equation}\label{commute}
    \begin{aligned}
        &u_{i\bj k}-u_{ik \bj}= -g^{l\bm} R_{k\bj i \bm}u_l,\\
        &u_{i\bj k}-u_{ki \bj}= -g^{l\bm} R_{i \bj k \bm}u_l,\\
        &u_{i\bj k\bl}-u_{k\bl i\bj}=g^{p\bq}(R_{k\bl i \bq}u_{p\bj}-R_{i\bj k \bq}u_{p\bl})
    \end{aligned}
\end{equation}

\noindent for $R_{i \bj k \bl}$ being the components of $(4,0)$ Riemannian curvature tensor.

We write the local coordinate formula for $\Omega + \sqrt{-1} \pbp u \wedge \omega^{p-1} $ as follows. Let $ \sqrt{-1}\partial \bpartial u \wedge \omega^{p-1} 
     = (p-1)! (\sqrt{-1})^{p^2} \sum\limits_{I, J \in \fI_p} U_{I\bJ} dz_I \wedge d \bz_J .$ Then in local orthonormal coordinates

$$ U_{I\bJ} =   \begin{cases}
    \sum\limits_{i \in I} u_{i\bi} 
          \;\; &I = J \\
       (-1)^{(i|I)+(j|J)} {u_{i\bj}}   &|I \cap J| = p-1, \; i \in I \setminus J, \; j \in J \setminus I \\
    0 & |I \cap J| \leq p-2.
    \end{cases}
    $$

Clearly in a general coordinate system, there will be more terms involving $u_{i \bj}$ and cofactors of order $(p-1)$ of $(g_{i \bj})$ in the above formula. We will not write the general form, as these extra terms do not play a significant role in the calculations. The only part this needs to be addressed is in \eqref{s5}, where we control the extra terms by  $C \lambda_1 \sum\limits_{I}Z^{I \bI}$.

We will use the notation $X_{I \bJ}$ for the matrix associated to the background $(p,p)$-form $p.\Omega$, so that

$$\Omega_{I \bJ} = (p-1)! \sqrt{-1}^{p^2} X_{I \bJ},$$

\noindent and

$$Z_{I \bJ} = X_{I \bJ} + U_{I \bJ}$$

\

\noindent denotes the matrix associated to the form $p\Omega_u$.

\subsection{\bf Concavity of the operator:} We show that the operator 

$$G(D^2 u) = \log{\det(X_{I \bJ} + U_{I \bJ})}$$ is concave as a function of $D^2 u$. All the calculations are done pointwise in an orthonormal coordinate system. Denote 

$$G^{i \bj} = \frac{\p G}{\p u_{j \bi}}, \hspace{.7 cm} G^{i \bj, k\bl} = \frac{\p^2 G}{\p u_{j \bi} \p u_{l\bk}}.$$

Then the first derivatives of $G$ has the following form.

$$G^{i \bj} = Z^{I \bJ} \frac{\p Z_{J \bI}}{\p u_{j \bi}} = \sum\limits_{\{I':\; i, j \notin I'\}} (-1)^{(i|I'_i)+(j|I'_j)} Z^{I'_i \bI'_j},$$

\noindent and

\begin{align}
G^{i \bj, k\bl} &= - \sum\limits_{\{I':\; i, j \notin I'\}} (-1)^{(i|I'_i)+(j|I'_j)} Z^{I'_i \bK }Z^{J \bI'_j} \frac{\p Z_{K \bJ}}{\p u_{l \bk}}\\
& = - \sum\limits_{\{I':\; i, j \notin I'\}} \sum\limits_{\{J':\; k, l \notin J'\}} (-1)^{(i|I'_i)+(j|I'_j) + (k|J'_k)+(l|J'_l) }Z^{I_i' \bJ'_l} Z^{J'_k \bI'_j}. \nonumber
\end{align}

For concavity, it is enough to show that this operator is negative-definite. Consider the matrix $(W^{I \bJ})$ that diagonalizes $(Z^{I \bJ})$.

$$Z^{K \bL} = W^{K\bM} \Lambda^{-1}_{M} W^{M \bL},$$

Then for any non-zero Hermitian matrix $\eta = (\eta_{i \bj})$,
\begin{equation}
\begin{aligned}
       G^{i \bj, k \bl} \eta_{j \bi} \eta_{l \bk} &=  - \sum\limits_{I', J'} \sum\limits_{i,j \notin I'} \sum\limits_{k,l \notin J'} (-1)^{(i|I'_i) + (k | J'_k)+ (j|I'_j) + (l| J'_l)}Z^{I_i' \bJ'_l}  Z^{J'_k \bI'_j} \eta_{j \bi} \eta_{l \bk}\\
       & = -\sum_{M,N} \Lambda_{M }^{-1} \Lambda_{N}^{-1}  \left| \sum_{I'}\sum_{\{i, j \notin I'\}} W^{I'_i \bM}(-1)^{(i|I'_i)}\eta_{j\bi}(-1)^{(j|I'_j)} W^{N \bI_j'}\right|^2 < 0
\end{aligned}    
\end{equation}


This would imply that $G$ is concave in $D^2 u$.

\





For ease of notation denote

$$\sigma(i,j,k,l) = (-1)^{(i|I'_i)+(j|I'_j) + (k|J'_k)+(l|J'_l)} \;\;\text{and }\;\;\sigma(i,j) = (-1)^{(i|I'_i)+(j|J'_j)}.$$

The multi-indices may differ in this notation, but should be clear from the context.

\




For the second order estimates, we could use the fact that $0$ is a $\cC$-subsolution to this equation as defined in \cite{Szekelyhidi18}. But since this equation has the form of a Monge-Amp\`ere equation, it is simpler to show that in orthonormal coordinates that diagonalizes the Hessian,

\begin{equation}\label{trace-Z}
\begin{aligned}
    N &= Z^{I \bJ} Z_{J \bI} =  Z^{I \bJ} X_{J \bI} + Z^{I \bJ} U_{J \bI}\\
    &=Z^{I \bJ} X_{J \bI} +Z^{I'_i \bI'_i} u_{i \bi},
\end{aligned}   
\end{equation}

\noindent so that

\begin{equation} \label{A-coefficient}
    - Z^{I'_i \bI'_i} u_{i \bi} = Z^{I \bJ} X_{J \bI}-N \geq \alpha \sum_I Z^{I \bI} -N
\end{equation}

\noindent where $\alpha$ is a positive constant, since $X$ is positive-definite. In the last line we use the basic inequality from linear algebra that for two symmetric positive-definite matrices $A$ and $B$ with eigenvalues $\lambda_1 \geq \hdots \geq \lambda_n$, and $\mu_1 \geq \hdots \geq \mu_n$, we have the following inequality.

$$\sum_{i, j} A^{i j} B_{i j} \geq \sum\limits_{i =1}^n \lambda_i \mu_{n-i}.$$

It should also be clear that $\sum Z^{I \bI} \geq \gamma>0$, for some positive constant $\gamma$ depending on the function $F$ in \eqref{pde1}.

\subsection{\bf Ellipticity of the operator $G$:}

We show that for any $u$ such that $\Lambda(\Omega_u) \in \Gamma_N$, the operator $G$ is elliptic. Consider any non-zero vector $\xi \in \mathbb C^n$. 

\begin{equation}\label{ellipticity}
    \begin{aligned}
        \sum\limits_{i,j } G^{i \bj} \xi_j \bar{\xi}_i &= \sum\limits_{i,j} \sum\limits_{\{I':\;i,j \notin I'\}}  (-1)^{(i|I'_i)+(j|I'_j)}Z^{I'_i \bI'_j} \xi_j \bar{\xi}_i \\ & = \sum_{I'} \sum\limits_{i, j \notin I'} Z^{I'_i \bI'_j} ((-1)^{(j|I'_j)}\xi_j) \ol{( (-1)^{(i|I'_i)}{\xi}_i)} \geq 0
    \end{aligned}
\end{equation}

\noindent since for each $I'$, the matrix $(Z^{I'_i \bI'_j})_{i , j}$ is an $(n-p+1) \times (n-p+ 1)$ principal submatrix of $(Z^{I \bJ})$, and is positive-definite.

\section{$C^0$ estimates}

The $C^0$ estimate can be obtained by using the ABP maximum principle. In \cite{Blocki11}, Blocki derived $C^0$ estimates for the complex Monge-Amp\`ere equation using the ABP maximum principle. This was later extended to more general fully nonlinear equations by Szekelyhidi \cite{Szekelyhidi18} on Hermitian manifolds. We say that a smooth function $\ul{u}$ is a $\cC$-subsolution for 
    $$f(\Lambda) = \log(\det(\Omega_u)) = F(x)$$
    if at each $x \in M$, the set $\{\Lambda( \Omega_{\ul{u}}) + \Gamma_N \} \cap \partial \Gamma^{F(x)}$ is bounded. Here 
    $$\partial \Gamma^{F(x)} = \{\Lambda\in \mathbb R^N: f(\Lambda) = F(x)  \}$$
\noindent is the level set at $F(x)$ of the function $f$. Note that the constant $b$ in\eqref{pde1} is uniformly bounded by the maximum principle.

From the definition, it can be seen that $0$ is a $\mathcal{C}$-subsolution to $\eqref{main-pde}$. A smooth function $\ul{u}$ being a $\mathcal C$-subsolution implies that there exist independent constants $\delta > 0$ and $R >0$ such that at each $x$

\begin{equation}\label{Csub}
    \{\Lambda( \Omega_{\ul{u}}) - \delta \mathbf{1} + \Gamma_N \} \cap \partial \Gamma^{F(x)}  \subset B_R(0)
\end{equation}

\

First observe that $\tr_{\omega} (Z_{I \bJ}) > 0$, would give the elliptic equation

\begin{equation}\label{trace}
 a^{i \bj} u_{i \bj} + g(z) > 0
\end{equation}




We normalize the solution such that $\sup\limits_M u =0$. So it is enough to estimate $\inf\limits_M u$. Let $ m = \inf\limits_M u$ be attained at a point $z_0$ on $M$ and assume that $m < 0$. Choose coordinates that takes $z_0$ to the origin and consider the coordinate ball $B(1) = \{z : |z| <1\}$. Let $v = u + \kappa |z|^2$ for a small $\kappa > 0$, so that $\inf\limits_{B(1)} v = m = v(0)$, and $\inf\limits_{z \in \partial B(1)} v(z) \geq m + \kappa$. Also 

Then by the ABP maximum principle for upper contact sets (\cite{Szekelyhidi18}, Proposition $10$), the set 

$$\Gamma^+_{\kappa} = \{x \in B(1): |Dv(x)| \leq \frac{\kappa}{2} \text{ and } v(y) - v(x) \geq Dv(x).(y-x) \text{ for all } y \in B(1)  \}$$

\noindent satisfies 

\begin{equation}\label{ABPint}
    \int_{\Gamma^+_{\kappa}} \det{D^2 v} \geq c_0 \kappa^{2n}. 
\end{equation}
\noindent for some positive constant $c_0$. The following can be verified at any $x \in \Gamma^+_{\kappa}$ as in \cite{Blocki11},

\begin{enumerate}[(i)]
    \item $D^2v(x) \geq 0$.
    \item \label{item2} $\det(D^2v) \leq 2^{2n} (\det v_{i \bj})^2$.
    \item $u_{i \bj} (x) \geq - \kappa \delta_{i j}$.
\end{enumerate}

It follows from $v = u+ \kappa |z|^2$ that at $x \in \Gamma_{\kappa}^+$

$$\Lambda(\Omega_{ u}) \in \{\Lambda( \Omega_{\ul{u}}) - \delta \mathbf{1} + \Gamma_N \}$$

\noindent for $\delta$ small depending on $\kappa$.

From the equation \eqref{pde1}, it is clear that $\Lambda(\Omega_{ u}) \in \partial \Gamma^\sigma$ for $\sigma = e^{F} \det \omega$. Hence it follows from \eqref{Csub} that $\Omega_{u}$ is bounded and 

\begin{equation}\label{wedge-bound}
    |\sqrt{-1}\pbp u \wedge \omega^{p-1}|_{\omega} \leq C.
\end{equation}

This would imply that $|u_{i \bj}| \leq C$. To see this, choose orthonormal coordinates such that $u_{i \bj} = \sum\limits_i\lambda_i \delta_{i j}$ and $g_{i \bj} = \delta_{i j}$ at a point. Then \eqref{wedge-bound} implies that $\sum\limits_{i \in I} \lambda_i$ is bounded for any set $I$ with $|I| = p$. This combined with $\sum \lambda_i + C > 0$ gives the required bound.

Note that this argument can be simplified in the Monge-Amp\'ere case \cite{Blocki11}, but this is more general. This gives a uniform bound for $v_{i \bj}$ in $\Gamma_{\kappa}^+$. Then from \eqref{item2} and \eqref{ABPint},

$$c_0 \kappa^{2n} \leq C' \mathrm{vol}(\Gamma_{\kappa}^+),$$

\noindent for some constant $C' >0$. From the weak Harnack inequality (\cite{GT83}, Theorem $8.18$) applied to equation \eqref{trace}, on $B(1)$

$$  \int_{B(1)} |u|^p \leq C (1 + \inf_{B(1)}|u|) \leq C'$$

This implies that $|v|_{L^p}$ is bounded. In the set $\Gamma_{\kappa}^+$, $v(x) \leq v(0) + \dfrac{\kappa}{2} < 0$. Putting these together with the following calculation

\begin{equation}
    \mathrm{vol}(\Gamma^+_{\kappa}) \left|v(0) + \frac{\kappa}{2} \right|^p ~\leq ~\int_{\Gamma^+_{\kappa}} |v|^p ~\leq ~C,
\end{equation}

\noindent we get $|m + \kappa/2|^p \leq C \kappa^{-2n}$, which shows that $m$ is bounded.

\ 

\

\section{Second-order estimates}

\

In this section, we will obtain estimates for the complex Hessian of $u$ that depends on the gradient $\nabla u$ in a quadratic way. Define $\lambda = (n-p+1) C_0 + \Delta_{\omega} u$, where $C_0 = 1 + \sup_M {|X|}$. Consider the following test function.

$$Q = \log(\lambda) + \varphi(|\nabla u|^2) + \psi(u).$$

\

 The function $\vp(t)$ is chosen so that $\vp'' = 2 (\vp')^2$. For example, similar to \cite{TW17}, we choose

$$\varphi(t) = -\frac{1}{2} \log{\left(1 - \frac{t}{2\sup{(1 + |
\nabla u|^2)}}\right)}.$$

This is well-defined for $t \in [0, \sup(1 + |\nabla u|^2)-1]$, and satisfies

\begin{equation}\label{s1}
\frac{1}{4 \sup(1 + |\nabla u|^2)}\leq \vp' \leq \frac{1}{2 \sup (1 + |\nabla u|^2)}.
\end{equation}

Let $\psi(t)$ defined on $[\inf u , 0]$ be

$$\psi(t) = -A \log{\left(1 + \frac{t}{2 (1 + \sup |u| )}\right)}$$

\noindent for a large constant $A>0$ to be defined later, so that 

\begin{equation}\label{s2}
{\psi'}^2 \leq \frac{A^2}{(2 + \sup |u|)^2 }, \hspace{0.5cm} \text{and} \hspace{0.5cm}\psi'' \geq \frac{A}{(2 + 2 \sup |u|)^2}.
\end{equation}

Assume that $Q$ attains maximum at some $z_0$ in $M$. Consider a local holomorphic coordinate system around $z_0$ such that at $z_0$, the complex Hessian of $u$ is diagonalized with eigenvalues $\lambda_1 \geq \hdots \geq \lambda_n$, and 

$$g_{i \bj}(z_0) = \delta_{i j}, \hspace{0.5cm} dg_{i \bj}(z_0) = 0,$$
\noindent for all $i,j$.

\begin{lemma}\label{comparison-lemma} If $\lambda \gg 1$, then

    $$\frac{\lambda}{2n} \leq\lambda_1 \leq \lambda.$$
\end{lemma}

\begin{proof}

Since $Z$ is positive-definite, we have that $Z_{I \bI} >0$ for each $I$. Sum this inequality over all the multi-indices $I$ such that $1 \notin I$, to get

\begin{equation}
    \begin{aligned}
        0 < \sum_{\{I|1 \notin I\}} X_{I \bI} + \sum_{\{I|1 \notin I\}} \sum\limits_{i \in I}\lambda_i& \leq {n-1\choose p}(C_0-1) + {n-2\choose p-1} \sum_{i \neq 1} \lambda_i.
    \end{aligned}
\end{equation}

This implies

\begin{equation}
\begin{aligned}
   \lambda &= \lambda_1 + \sum\limits_{i \neq 1} \lambda_i + (n-p+1) C_0  \\ 
   &\geq \lambda_1 - (C_0-1) {n-1\choose p}{n-2\choose p-1}^{-1} + (n-p+1)C_0\\
   &=\lambda_1 + \frac{(n-p)(p-1) + 1}{p}C_0 + \frac{(n-1)}{p}\geq \lambda_1
   \end{aligned}
\end{equation}

For the other side, we have

\begin{equation}
    \begin{aligned}
        \lambda &\leq n \lambda_{1} + (n-p+1) C_0\\
    \end{aligned}
\end{equation}

Now assuming $\lambda > 2 (n-p+1) C_0$ gives $\lambda_1 \geq \dfrac{\lambda}{2n}$

\end{proof}


All calculations are performed in this coordinate system and uses the Levi-Civita connection associated to the metric $\omega$. We will show that at the maximum point $z_0$ of $Q$, $\lambda$ is bounded by $C \sup (1 + |\nabla u|^2)$, for some independent constant $C$. This will be enough to show that $|\sqrt{-1} \pbp u|$ is uniformly bounded. Indeed if $\lambda_1 \leq  \lambda \leq C \sup (1 + |\nabla u|^2)$, then positive-definiteness of $(Z_{I \bJ})$ implies a lower-bound of the same form for the smallest eigenvalue, $|\lambda_n| \leq C \sup{(1 +|\nabla u|^2)}$. 

For contradiction, assume that $\lambda \gg \sup (1 + |\nabla u|^2)$. In the calculations, we will implicitly use the fact that $\lambda$ can be assumed to be larger than any a priori bounded constant.






We split the indices in $\fI_p$ into two types. Let

$$\cP = \{ I \in \fI_p : \; (1+ 2\delta)\lambda\geq |Z_{I \bI}| \geq \delta\lambda\}.$$

\noindent where the constant $\delta>0$ will be chosen later. Also let

$$\cP'= \{I \in \fI_p : \;  |Z_{I \bI}| <\delta \lambda\}$$

Assume first that $\cP \neq \emptyset$, and $\cP' \neq \emptyset$. Define the set

$$S = \{k: |\lambda_k| \geq \tau \lambda\}$$

\noindent where $1 > \tau>0$ is a small constant to be fixed later that depends on $\delta$, $A$ and other known quantities, and $\tau< \dfrac{1}{2n}$. Clearly $1 \in S$.

\



We calculate $ \nabla_a \nabla_{\bar a}$ of $\log{\det Z}$.

\begin{equation}\label{s4}
\begin{aligned}
    \nabla_{\ba} \log{\det Z} &= Z^{I \bJ} \nabla_{\bar a} Z_{J \bI}\\
    &= Z^{I \bJ} \frac{\partial U_{J \bI}}{\partial u_{l \bk}} u_{l \bk \ba} + Z^{I \bJ} \nabla_{\ba} X_{J \bI}\\
    & = \sum\limits_{k,l}\sum\limits_{\{I': \;  l,k \notin I'\}}  \sigma(k,l) Z^{I'_k \bI'_l} u_{l \bk \ba} + Z^{I \bJ} \nabla_{\ba} X_{J \bI}
\end{aligned}
\end{equation}

In the following calculations $\sum\limits_{k,l}\sum\limits_{\{I': \;  l,k \notin I'\}}$ will be omitted when clear.

\begin{equation}\label{s5}
\begin{aligned}
    \nabla_a\nabla_{\bar a} \log{\det Z} & \leq Z^{I \bJ} \nabla_a\nabla_{\bar a} X_{J \bI} - Z^{I \bJ} Z^{K \bL} \nabla_{a} Z_{L \bI} \nabla_{\ba} Z_{J \bK} +  \sigma(i,j) Z^{I'_i \bI'_j}u_{j \bi a \ba} + C\lambda\sum\limits_{I}Z^{I \bI}.\\
\end{aligned}
\end{equation}


 Applying the covariant derivative $\nabla_{\ba}$ to equation \eqref{main-pde}

\begin{equation}\label{s4.1}
\begin{aligned}
        \nabla_{\ba} F = \sigma(i,j) Z^{I'_i \bI'_j} u_{j \bi \ba} + Z^{I \bJ} \nabla_{\ba} X_{J \bI}
\end{aligned}
\end{equation}

Applying the maximum principle to the function $Q$ at $z_0$,

\begin{equation}\label{s6}
\begin{aligned}
    0 &\geq  \sigma(i,j) Z^{I'_i \bI'_j} \sum\limits_a\frac{u_{a \ba j \bi}}{\lambda}- \sigma(i,j) Z^{I'_i \bI'_j}\frac{\nabla_j\lambda \nabla_{\bi}\lambda}{\lambda^2} +  2 \vp '    \fRe{ \left\{\sigma(i,j)Z^{I'_i \bI'_j} u_{\bk} u_{k j \bi}\right\}}  + \mathcal{H} 
\end{aligned}
\end{equation}

\noindent with

\begin{equation}\label{s7}
\begin{aligned}
    \mathcal{H} &= \sigma(i,j) Z^{I'_i \bI'_j} \varphi''   w_{j} w_{ \bi} +  \sigma(i,j) Z^{I'_i \bI'_j}\vp'\sum_k(u_{\bk j}u_{k \bi} + u_{k j}u_{\bk \bi} )\\
    &+ \psi'  \sigma(i,j) Z^{I'_i \bI'_j}u_{j \bi} + \psi'' \sigma(i,j) Z^{I'_i \bI'_j} u_j u_{\bi}
\end{aligned}
\end{equation}

\noindent where we set

$$w_{j} = \sum_k(u_{\bk} u_{k j} + u_{k} u_{\bk j})$$
$$w_{ \bi} = \sum_k(u_{\bk} u_{k \bi} + u_{k} u_{\bk \bi})$$

By using $\nabla Q = 0$ at $z_0$, we get the following compatibility equations.

\begin{equation}\label{s8}
    \begin{aligned}
    &\frac{\nabla_j\lambda}{\lambda} + \vp'w_{ j} + \psi' u_j = 0 \\
    & \frac{\nabla_{\bi}\lambda}{\lambda} + \vp'w_{ \bi} + \psi' u_{\bi} = 0
    \end{aligned}
\end{equation}

We commute the derivatives in the fourth order terms and control the resulting curvature terms as follows.

\begin{equation} \label{s9}
    u_{a \ba j \bi} = u_{j \bi a \ba} + \sum_p(R_{j \bi a \bp}u_{p \ba} -  R_{a \ba j \bp} u_{p \bi}) \geq u_{j \bi a \ba} - C \lambda
\end{equation}

For the third-order terms

\begin{equation} \label{s10}
    u_{j \bi k} = u_{k j \bi} - \sum_l R_{j \bi k \bl} u_{l} .
\end{equation}

Combining \eqref{s4.1}, \eqref{s6}, \eqref{s9}, and \eqref{s10} gives

\begin{equation}\label{s11}
\begin{aligned}
0 \geq &  \sigma(i,j) Z^{I'_i \bI'_j} \sum\limits_a\frac{u_{j \bi a \ba}}{\lambda}- \sigma(i,j) Z^{I'_i \bI'_j} \frac{\nabla_j\lambda \nabla_{\bi}\lambda }{\lambda^2} +  \mathcal{H} - C \vp' (1 + |\nabla u|^2) \sum Z^{I \bI} - C( 1+ \sum Z^{I \bI})
\end{aligned}
\end{equation}





\










Combining equations \eqref{s11} and \eqref{s5} gives

\begin{equation} \label{s14}
\begin{aligned}
    0 & \geq  \frac{1 }{\lambda} \sum_a Z^{I \bJ} Z^{K \bL} \nabla_{a} Z_{L \bI} \nabla_{\ba} Z_{J \bK}- \sigma(i,j) Z^{I'_i \bI'_j}\frac{\nabla_{j}\lambda \nabla_{\bi}\lambda}{\lambda^2} +  \mathcal{H} - C (1+ \sum Z^{I \bI}) \\
    &  - C \vp' (1 + |\nabla u|^2) \sum Z^{I \bI}
\end{aligned}
\end{equation}

The following lemmas play an important role in the calculation.

\begin{lemma} \label{inversion-lemma}
    Given a positive-definite $n \times n$ Hermitian matrix $A = (a_{i \bj})$ with eigenvalues $\lambda_1, \hdots, \lambda_n$, and inverse $(a^{i \bj})$ the following are true. For any $i$

\begin{equation} \label{s15}
   a^{i \bi} \geq \frac{1}{a_{i \bi}}.
\end{equation}

 For any $\epsilon>0$ small, we have the upper bound

\begin{equation}\label{s16}
a^{i \bi} \leq \frac{1+ \epsilon}{a_{i \bi}} + \frac{C}{a_{i \bi}^2}\sum_{j \neq i} a^{j \bj},
\end{equation}

\noindent for a constant $C$ that depends only on $\epsilon$ and the non-diagonal entries of $A$.

\end{lemma}

\begin{proof}
    Given orthonormal eigenvectors $V_i$ corresponding to each eigenvalue $\lambda_i$, we can write
    \begin{align}
    a_{i \bi} = \sum_k |V_{i}^k|^2 \lambda_k,\hspace{0.5cm} \text{and} \hspace{0.5cm} a^{i \bi}  = \sum_j |V_{i}^j|^2\frac{1}{\lambda_j}.
    \end{align}

Now \eqref{s15} follows from

\begin{equation} \label{s17}
    a_{i \bi} a^{i \bi} = \sum\limits_{k ,j } \frac{\lambda_k |V_i^k|^2 |V_i^j|^2}{\lambda_j} \geq 1.
\end{equation}

For the second inequality, from $A^{-1} A = I$,

$$a^{i \bi }a_{i \bi} = 1- \sum\limits_{j \neq i}a^{i \bj} a_{j \bi} \leq 1 + \epsilon a_{i \bi} a^{i \bi } + C(\epsilon) \sum\limits_{j \neq i}\frac{a^{j \bj}}{a_{i \bi}}$$

  Here we used $|a^{i \bj}| \leq \sqrt{a^{i \bi} a^{j \bj}}$.  The result follows.
\end{proof}

The lower bound \eqref{s15} is not sufficient. So we prove some stronger lemmas that hold in the context of this problem. Let $V_{I_1}, V_{I_2}, \hdots, V_{I_N}$ be orthonormal eigenvectors corresponding to the eigenvalues $\Lambda_{I_1}, \hdots, \Lambda_{I_N}$ of the matrix $(Z_{I \bJ})$. Recall the definition $\cP = \{I: \delta \lambda \leq |Z_{I \bI}| \leq (1 + 2\delta) \lambda\}$. Let $\delta' $ be larger than $2\delta$. Then define

$$\cT = \{K:\; \Lambda_K \geq (1 + \delta') \lambda \}.$$

\

\begin{lemma}\label{evlemma1}\begin{enumerate}
   \item \label{evlemma1.1}  For all $K \in \cT$, and $I$ such that $Z_{I \bI} \leq (1+ 2 \delta) \lambda$

   \ 

\begin{equation}
    |V_K^I| \leq \frac{c}{\lambda}
\end{equation}

\

\noindent for a constant $c$ depending on $\delta' - 2\delta$. More precisely $c = O((\delta' - 2 \delta)^{-1})$.

\

    \item \label{evlemma1.2}   For $I$ or $J$ as above, and $I \neq J$

    \ 
    
          \begin{equation} \label{evlemma_bound}
              \begin{aligned}
                  &|\sum_{K \notin \cT} V_K^I \ol{V_K^J}| \leq \epsilon_3\\
                  &\sum_{K \notin \cT} |V_K^I|^2 \geq (1 - \epsilon_3)
              \end{aligned}
          \end{equation}

          \

\noindent for a constant $\epsilon_3 \leq \dfrac{c}{\lambda}$.   
\end{enumerate}
\end{lemma}

\

\begin{proof}

For part \eqref{evlemma1.1}, it follows from $ Z. V_K = \Lambda_K V_K$ that
  \begin{equation} \label{s18}
        |\sum\limits_{J \neq I}Z_{I \bJ} V^J_K |\geq \Lambda_K |V_K^I| - Z_{I \bI} |V^I_K| \geq ( \delta' - 2\delta) \lambda |V_K^I|
  \end{equation}
  
This gives $|V_K^I| \leq \dfrac{c}{ \lambda}$, for $c$ depending only on $(X_{I \bJ})$ and $\delta' - 2\delta$.

\

    Part \eqref{evlemma1.2} now follows from part \eqref{evlemma1.1} since $\sum_K V^I_K \ol{V_K^J} = \delta_{I J} .$

\end{proof}




We may assume that $\delta' = 3 \delta$. The constant $\delta$ will be chosen in \eqref{fixalpha} and \eqref{fixconstants} to depend on $\sup |u|$ and the constant $A$, both of which are independently known.




    










\

\

\

\subsection{The term $  \sum\limits_a Z^{I \bJ} Z^{K \bL} \nabla_{a} Z_{L \bI} \nabla_{\ba} Z_{J \bK}$}

\ 

\

In this section, we obtain a comparison between the concavity terms $ \sum_a Z^{I \bJ} Z^{K \bL} \nabla_{a} Z_{L \bI} \nabla_{\ba} Z_{J \bK}$ and the term $Z^{I'_i \bI'_j} \dfrac{\nabla_j \lambda \nabla_{\bi} \lambda}{\lambda}$. Define the complex vector $\xi_I \in \mathbb C^N$ with components

$$\xi_I^{M,a} = \frac{1}{\sqrt{\Lambda_M}} \sum\limits_J V_M^J \nabla_{\ba }Z_{I \bJ}.$$

By spectral decomposition

\begin{equation}\label{spectral-decomp}
    \begin{aligned}
        \sum_a Z^{I \bJ} Z^{K \bL} \nabla_{a} Z_{L \bI} \nabla_{\ba} Z_{J \bK} & = \sum_a\sum\limits_{M,K}\frac{\left| V_M^I \nabla_a Z_{L \bI} \ol{V_K^L} \right|^2}{\Lambda_M \Lambda_K}\\
        & = \sum_a\sum\limits_{ M,K } \frac{\left|\sum\limits_{I} V_K^{I}\xi_{I}^{M,a} \right|^2}{\Lambda_K}.
    \end{aligned}
\end{equation}

We define the pairing $(I',i)$ when $i \notin I'$. Then let

$$\mathfrak{J} = \{(I',i) \in \fI_{p-1}\times \{1, \cdots, n\}: i \in S^c, I'_i \in \cP', 1 \notin I'\; \}.$$

If $(I',i) \in \mathfrak J$, then

$$\sum\limits_{j\in I'} \lambda_j \leq (\delta + \tau)\lambda + C_0.$$

\begin{lemma}\label{lemma-concavity}
  \begin{equation}
      \begin{aligned}
          \sum_{a}\sum\limits_{M,K}\frac{\left| V_M^I \nabla_a Z_{L \bI} \ol{V_K^L} \right|^2}{\Lambda_M \Lambda_K} \geq &\frac{1}{(1+ \delta_0) \lambda} \sum\limits_M \sum_{I' \in \mathfrak{K}} \left| \sum\limits_{a \notin I'} (-1)^{(a| I'_a)} \xi_{I'_a}^{M,a}\right|^2 \\
          &+ (1-\epsilon')\sum_M\sum_a\sum\limits_{K \in \cT} \frac{\left|\sum\limits_{I} V_K^{I}\xi_{I}^{M,a} \right|^2}{\Lambda_K}.
      \end{aligned}
  \end{equation}

  \noindent where $\delta_0 =  (n-p)(\delta + \tau) + \dfrac{c}{\epsilon'\lambda^2}$, and the set $\mathfrak{K} = \{I': (I',i) \in \mathfrak J \text{ for some } i \}$. Here $\epsilon'>0$ is a small constant that can be chosen to be $\lambda^{-1}$.
\end{lemma}

\begin{proof}

\begin{equation}\label{T-split}
    \begin{aligned}
       &\left(\sum_{I' \in \mathfrak{K}} \left| \sum\limits_{a \notin I'} (-1)^{(a| I'_a)} \xi_{I'_a}^{M,a}\right|^2\right)^{1/2} =  \left(\sum_{I' \in \mathfrak{K}} \left| \sum\limits_{a \notin I'} (-1)^{(a| I'_a)} \sum\limits_{I,K} V_K^{I'_a}V_K^{I}\xi_{I}^{M,a}\right|^2\right)^{1/2}\\
       & =   \left(\sum_{I' \in \mathfrak{K}} \left| \sum\limits_{a \notin I'} \sum\limits_{K \in \cT}(-1)^{(a| I'_a)} \sum\limits_{I} V_K^{I'_a}V_K^{I}\xi_{I}^{M,a} + \sum\limits_{a \notin I'} \sum\limits_{K \notin \cT}(-1)^{(a| I'_a)} \sum\limits_{I} V_K^{I'_a}V_K^{I}\xi_{I}^{M,a}\right|^2 \right)^{1/2}\\
       & \leq \left(\sum_{I' \in \mathfrak{K}} \left| \sum\limits_{a \notin I'} \sum\limits_{K \in \cT}(-1)^{(a| I'_a)} \sum\limits_{I} V_K^{I'_a}V_K^{I}\xi_{I}^{M,a} \right |^2\right)^{1/2}\\
       &\;\;\;+ \left(\sum_{I' \in \mathfrak{K}}\left|\sum\limits_{a \notin I'} \sum\limits_{K \notin \cT}(-1)^{(a| I'_a)} \sum\limits_{I} V_K^{I'_a}V_K^{I}\xi_{I}^{M,a}\right|^{2}\right)^{1/2}\\
    \end{aligned}
\end{equation}

\noindent by the triangle inequality. To the first term in the last sum, we apply the Cauchy-Schwarz inequality for indices $K$ and $a\notin I'$ with weights $\Lambda_K$.

\begin{equation}\label{cauchy1}
    \begin{aligned}
      \sum\limits_{I' \in \mathfrak K}\left|\sum\limits_{K \in \cT, a \notin I'} \left((-1)^{(a| I'_a)} V_K^{I'_a}\right)\right. & \left.\left(\sum\limits_{I} V_K^{I}\xi_{I}^{M,a} \right) \right|^2  \leq \sum\limits_{I' \in \mathfrak K}\left(\sum\limits_{K \in \cT} \Lambda_K\sum\limits_{a \notin I'}\left| V_K^{I'_a}\right|^2 \right)\\
      &\times\left( \sum\limits_{a \notin I',K \in \cT} \frac{\left|\sum\limits_{I} V_K^{I}\xi_{I}^{M,a} \right|^2}{\Lambda_K}\right)\\
       &\leq \left(  \sum\limits_{I' \in \mathfrak K}\sum\limits_{K \in \cT}  \sum\limits_{a \notin I'} \Lambda_K |V_K^{I'_a}|^2 \right)\left( \sum_a\sum\limits_{K \in \cT} \frac{\left|\sum\limits_{I} V_K^{I}\xi_{I}^{M,a} \right|^2}{\Lambda_K}\right)\\
       &\leq \frac{c}{\lambda}  \left( \sum_a\sum\limits_{K \in \cT} \frac{\left|\sum\limits_{I} V_K^{I}\xi_{I}^{M,a} \right|^2}{\Lambda_K}\right)
    \end{aligned}
\end{equation}

\noindent using \eqref{ineq2}. In the second line, we expanded the index set of $a$ to the set $\{1, \cdots,n\}$. This allows us to make the summation over $I' \in \mathfrak K$ to appear only in the first term. For $I' \in \mathfrak K$, and $a \notin I'$, we have

\begin{equation}
    \begin{aligned}
        Z_{I'_a \bI'_a} &= X_{I'_a \bI'_a} + \sum\limits_{j \in I'}\lambda_j + \lambda_a\\
        &\leq C_0-1+ (\delta + \tau) \lambda + \lambda_1 \leq (1+ 2 \delta) \lambda,
    \end{aligned}
\end{equation}

\noindent after setting $\tau < \frac{\delta}{2}$. Then by Lemma \ref{evlemma1},

\begin{equation}\label{ineq2}
    \begin{aligned}
     \sum\limits_{K \in \cT} \sum\limits_{I' \in \mathfrak K}  \sum\limits_{a\notin I'}\Lambda_K |V_K^{I'_a}|^2  \leq \frac{c}{\lambda}.
    \end{aligned}
\end{equation}

Consider the following matrix with columns indexed by $(a,K)$ for $K \notin \cT$ and rows by $I'\in \mathfrak K$, 

$$B_{I',(a,K)} = \begin{cases}(-1)^{(a| I'_a)} \sqrt{\Lambda_K}V_K^{I'_a}&, \;\; {a \notin I'}\\
0&, \;\; {a \in I'}
\end{cases}$$

\noindent and the vector $\chi_{(a,K)} = \dfrac{\sum\limits_{I} V_K^{I}\xi_{I}^{M,a} }{\sqrt{\Lambda_K}}$ for $K \notin \cT$ and a fixed $M$. Note that 

$$BB^* =  \Bigg[\sum\limits_{a \notin I' \cup J'}(-1)^{(a| J'_a) + (a| I'_a)}\sum\limits_{K \notin \cT} \Lambda_K V_K^{I'_a}\ol{ V_{K}^{J'_a}} \Bigg]_{I',J' \in \mathfrak K},$$

\noindent so that $|B|^2 \leq (1 + (n-p)(\delta + \tau))\lambda$. Then 

\begin{equation}\label{cauchy2}
    \begin{aligned}
      \sum\limits_{I' \in \mathfrak K}\left|\sum\limits_{K \notin \cT} \sum\limits_{a \notin I'}\left((-1)^{(a| I'_a)} V_K^{I'_a}\right)\right. & \left.\left(\sum\limits_{I} V_K^{I}\xi_{I}^{M,a} \right) \right|^2= |B\chi|^2 \\ 
      & \leq |B|^2 |\chi|^2\\
       & \leq (1 + (n-p)(\delta + \tau)) \lambda \left( \sum_a\sum\limits_{K \notin \cT} \frac{\left|\sum\limits_{I} V_K^{I}\xi_{I}^{M,a} \right|^2}{\Lambda_K}\right)
    \end{aligned}
\end{equation}

\noindent using the matrix inequality \eqref{matrix-ineq}. Notice that

\begin{equation}\label{matrix-ineq1}
    \begin{aligned}
        \sum\limits_{K \notin \cT} \Lambda_K V_k^{I'_a} \ol{V_K^{J'_a}} \preceq Z_{I'_a \bJ'_a} - \sum\limits_{K \in \cT} \Lambda_K V_K^{I'_a} \ol{V_K^{J'_a}} \preceq Z_{I'_a \bJ'_a}
    \end{aligned}
\end{equation}

For $I' \in \mathfrak K$, $(I',i) \in \mathfrak J$ for some $i$, which gives

\begin{equation}\label{matrix-ineq2}
    \begin{aligned}
        \sum\limits_{j \in I'} \lambda_j &= Z_{I'_i \bI'_i} - X_{I'_i \bI'_i} - \lambda_i\\
        & \leq (\delta+ \tau) \lambda.
    \end{aligned}
\end{equation}

Since $X_{I'_a \bJ'_a} \preceq C_0 \delta_{I' J'}$, using \cref{matrix-ineq1,matrix-ineq2} we have the following matrix inequality for submatrix indexed by $I', J' \in \mathfrak K$.

\begin{equation}\label{matrix-ineq}
    \begin{aligned}
         \Bigg[\sum\limits_{a \notin I' \cup J'}(-1)^{(a| J'_a) + (a| I'_a)}\sum\limits_{K \notin \cT} &\Lambda_K V_K^{I'_a}\ol{ V_{K}^{J'_a}} \Bigg]_{I',J' \in \mathfrak K}
          \preceq \left[\sum\limits_{a \notin I' \cup J'}(-1)^{(a| J'_a) + (a| I'_a)}Z_{I'_a \bJ'_a}\right]_{I',J' \in \mathfrak K}\\
         &= \left[\sum\limits_{a \notin I' \cup J'}(-1)^{(a| J'_a) + (a| I'_a)}X_{I'_a \bJ'_a} + \delta_{I'J'}\sum\limits_{a \notin I'}(\sum\limits_{i\in I'_a} \lambda_i)\right]_{I',J' \in \mathfrak K}\\
         &\preceq \left[\left((n-p+1)C_0 + \Delta u + (n-p)(\sum\limits_{i \in I'} \lambda_i)\right)\delta_{I'J'}\right]_{I',J' \in \mathfrak K}\\
         & \preceq (1 + (n-p)(\delta + \tau))\lambda \; \delta_{I'J'}.
    \end{aligned}
\end{equation}

Inserting \eqref{cauchy1} and \eqref{cauchy2} in \eqref{T-split}, we get

\begin{equation}
    \begin{aligned}
        &\left(\sum_{I' \in \mathfrak{K}} \left| \sum\limits_{a \notin I'} (-1)^{(a| I'_a)} \xi_{I'_a}^{M,a}\right|^2\right)^{1/2}  \leq \sqrt{(1 + (n-p)(\delta + \tau)) \lambda} \left( \sum_a\sum\limits_{K \notin \cT} \frac{\left|\sum\limits_{I} V_K^{I}\xi_{I}^{M,a} \right|^2}{\Lambda_K}\right)^{1/2}\\
        & \;\; + \sqrt{\frac{c}{\lambda}}  \left( \sum_a\sum\limits_{K \in \cT} \frac{\left|\sum\limits_{I} V_K^{I}\xi_{I}^{M,a} \right|^2}{\Lambda_K}\right)^{1/2}\\
        & \leq \sqrt{(1 + \delta_0) \lambda }\left[ \left( \sum_a\sum\limits_{K \notin \cT} \frac{\left|\sum\limits_{I} V_K^{I}\xi_{I}^{M,a} \right|^2}{\Lambda_K}\right) + \epsilon'\left( \sum_a\sum\limits_{K \in \cT} \frac{\left|\sum\limits_{I} V_K^{I}\xi_{I}^{M,a} \right|^2}{\Lambda_K}\right)\right]^{1/2}\\
        & \leq \sqrt{(1 + \delta_0) \lambda }\left[ \left( \sum_a\sum\limits_{K } \frac{\left|\sum\limits_{I} V_K^{I}\xi_{I}^{M,a} \right|^2}{\Lambda_K}\right) - (1- \epsilon')\left( \sum_a\sum\limits_{K \in \cT} \frac{\left|\sum\limits_{I} V_K^{I}\xi_{I}^{M,a} \right|^2}{\Lambda_K}\right)\right]^{1/2}
    \end{aligned}
\end{equation}

\noindent where we used the basic inequality $(\sqrt{a}x + \sqrt{b}y)^2 \leq \left(a + \dfrac{b}{\epsilon'}\right)(x^2 + \epsilon' y^2)$.

It follows this and \eqref{spectral-decomp} that 

\begin{equation}
    \begin{aligned}
        \sum_{a}\sum\limits_{M,N}\frac{\left| V_M^I \nabla_a Z_{L \bI} \ol{V_N^L} \right|^2}{\Lambda_M \Lambda_N} &=  \sum_M\sum_a\sum\limits_{ K } \frac{\left|\sum\limits_{I} V_K^{I}\xi_{I}^{M,a} \right|^2}{\Lambda_K}\\
        & \geq \frac{1}{(1+ \delta_0)\lambda} \sum\limits_M \sum_{I' \in \mathfrak{K}} \left| \sum\limits_{a \notin I'} (-1)^{(a| I'_a)} \xi_{I'_a}^{M,a}\right|^2\\
        &\;\;+ (1-\epsilon')\sum_M\sum_a\sum\limits_{K \in \cT} \frac{\left|\sum\limits_{I} V_K^{I}\xi_{I}^{M,a} \right|^2}{\Lambda_K}
    \end{aligned}
\end{equation}

This completes the proof. 
    
\end{proof}

We also have

\begin{equation}
    \begin{aligned}
 \sum\limits_{a \notin I'}(-1)^{(a| I'_a)}\xi_{I'_a}^{M,a} = \frac{1}{\sqrt{\Lambda_M}} \sum\limits_{a \notin I'}(-1)^{(a| I'_a)}\sum_{J} V_{M}^J \nabla_{\ba}(U_{I'_a \bJ} + X_{I'_a \bJ}).
    \end{aligned}
\end{equation}

We calculate

\begin{equation}
    \begin{aligned}
          \sum\limits_{a \notin I'}(-1)^{(a| I'_a)}\sum_{J} V_{M}^J \nabla_{\ba}U_{I'_a \bJ} &= \sum\limits_{a \notin I'}(-1)^{(a| I'_a)} \left(\sum\limits_{b \notin I', b\neq a} V_M^{I'_b} \nabla_{\ba} U_{I'_a \bI'_b} + V_{M}^{I'_a} \nabla_{\ba} U_{I'_a \bI'_a}\right.\\
          &\left.\;\;+ \sum_{m \in I', n \notin I', n \neq a}\sum\limits_{J = (I'_a\setminus m) \cup n}V_{M}^J \nabla_{\ba}U_{I'_a \bJ}\right) \\
          & = \sum\limits_{a \notin I'}\sum\limits_{b \notin I', b\neq a} (-1)^{(b| I'_b)}V_{M}^{I'_b}  \nabla_{\ba}u_{a \bar b} + \sum\limits_{a \notin I'}(-1)^{(a| I'_a)}V_{M}^{I'_a}  \sum\limits_{k\in I'_a}\nabla_{\ba} u_{k \bk} \\
          & =  \sum\limits_{a \notin I'}\sum\limits_{b \notin I', b\neq a} (-1)^{(a| I'_a)}V_{M}^{I'_a}  \nabla_{\bar b}u_{b \bar a} +\sum\limits_{a \notin I'} (-1)^{(a| I'_a)}V_{M}^{I'_a}  \sum\limits_{k\in I'_a}\nabla_{\ba} u_{k \bk}\\
          & = \sum\limits_{a \notin I'} (-1)^{(a| I'_a)}V_{M}^{I'_a} \nabla_{\ba} \lambda 
    \end{aligned}
\end{equation}

\noindent where we use the commutation formula $\nabla_{\bar b}u_{b \ba} = \nabla_{\bar a} u_{b \bar b}$ on K\"ahler manifolds and interchanged $a$ and $b$ indices in the third line.

Here the term

\begin{equation}
    \begin{aligned}
       \sum\limits_{a \notin I'} (-1)^{(a| I'_a)} \sum_{\substack{m \in I',\\ n \notin I',\\ n\neq a}}\sum\limits_{J = (I'_a\setminus m) \cup n}V_{M}^J \nabla_{\ba}U_{I'_a \bJ}
        & =  \sum\limits_{m \in I'}\sum\limits_{\substack{a,n \notin I',\\ a \neq n}}\sum\limits_{J = (I'_a\setminus m)\cup n} (-1)^{(m|I'_a)+ (a|I'_a)+ (n|J)}V^{J}_M \nabla_{\ba} u_{m\bn}\\
        & = \sum\limits_{m \in I'}\sum\limits_{a,n \notin I', a< n}\;\;\sum\limits_{J = (I'_a\setminus m)\cup n}[ (-1)^{(m|I'_a)+ (a|I'_a)+ (n|J)} \\
        &\;\;+ (-1)^{(m|I'_n)+ (n|I'_n)+ (a|J)}  ]V^{J}_M \nabla_{\ba} u_{m\bn}\\
        & = 0.
    \end{aligned}
\end{equation}

\noindent where in the second line we used that $\nabla_{\ba} u_{m\bn} =\nabla_{\bn} u_{m\ba}$ on K\"ahler manifolds and $J = (I'_a\setminus m) \cup n$ is symmetric in $a$ and $n$. In the last line we used that

$$(m|I'_a)+ (a|I'_a)+ (n|J) \equiv (m|I'_n)+ (n|I'_n)+ (a|J) + 1 \; (\text{mod } 2).$$

This can be seen by considering the three cases:

\begin{enumerate}
    \item If $m<a<n$,
$$(n|J)= (n|I'_n), \;(m|I'_a)=(m|I'_n), \; (a|I'_a) - 1=  (a|J).$$

    \item If $a<m<n$

    $$(n|J)= (n|I'_n), \; (m|I'_a)=(m|I'_n) + 1,\; (a|I'_a) =  (a|J).$$

    \item If $a<n<m$

    $$(n|J) - 1= (n|I'_n), \; (m|I'_a)=(m|I'_n), \; (a|I'_a)=  (a|J)$$
\end{enumerate}

\begin{equation}
    \begin{aligned}
        \sum_M\left|\frac{1}{\sqrt{\Lambda_M}} \sum\limits_{a \notin I'}(-1)^{(a| I'_a)} \sum\limits_J V_M^J \nabla_{\ba} X_{I'_a \bJ} \right|^2 \leq \sum_{M}\frac{C}{\Lambda_M} = C\sum_{I} Z^{I \bI}
    \end{aligned}
\end{equation}

It follows from the above that for $I' \in \mathfrak K$,

\begin{equation}
    \begin{aligned}
        \sum_M\left| \sum\limits_{a \notin I'} (-1)^{(a| I'_a)} \xi_{I'_a}^{M,a}\right|^2 &\geq (1-\epsilon_0)\sum_M\sum\limits_{a,b \notin I'} (-1)^{(a| I'_a)+  (b|I'_b)}\frac{V_M^{I'_a}\ol{V_M^{I'_b}}}{\Lambda_M}\nabla_b \lambda \nabla_{\bar a} \lambda\\
        &- \epsilon_0^{-1} \sum_M\left|\frac{1}{\sqrt{\Lambda_M}} \sum\limits_{a \notin I'} (-1)^{(a| I'_a)}\sum\limits_J V_M^J \nabla_{\ba} X_{I'_a \bJ} \right|^2\\
        &\geq (1-\epsilon_0) \sum\limits_{a,b \notin I'} \sigma(a,b)Z^{I'_a \bI'_{ b}}\nabla_b \lambda \nabla_{\bar a}  \lambda - C\epsilon_0^{-1} \sum_{I} Z^{I \bI}
    \end{aligned}
\end{equation}

Now it follows from \eqref{spectral-decomp} and Lemma \ref{lemma-concavity} that

\begin{equation}\label{third-order-control2}
    \begin{aligned}
      \sum_a Z^{I \bJ} Z^{K \bL} \nabla_{a} Z_{L \bI} \nabla_{\ba} Z_{J \bK} \geq \frac{(1-\epsilon_0)}{(1+ \delta_0) \lambda} \sum\limits_{I' \in \mathfrak K} \sum\limits_{a,b \notin I'} \sigma(a,b) Z^{I'_a \bI'_{ b}}\nabla_b \lambda \nabla_{\bar a}  \lambda - \frac{C}{\epsilon_0 \lambda} \sum_{I} Z^{I \bI}
    \end{aligned}
\end{equation}

\

\subsection{The term  $\dfrac{1}{\lambda^2} \sigma(b,l)Z^{I'_b \bI'_l}\nabla_l\lambda\nabla_{\bb}\lambda $}

 \ 
 
\



\

 Observe that when $(I',i)\in \fJ$

\begin{equation} \label{s20}
    Z_{I'_1 \bI'_1} = Z_{I'_i \bI'_i} + X_{I'_1 \bI'_1} - X_{I'_i \bI'_i} + \lambda_1 - \lambda_i.
\end{equation}

\

\noindent Lemma \ref{comparison-lemma} implies that

\begin{equation} \label{s21}
    \delta \lambda \leq Z_{I'_1 \bI'_1} \leq (1+ 2 \delta) \lambda
\end{equation}

\noindent assuming $0<\tau < \dfrac{\delta}{2}, \delta + \tau < \dfrac{1}{2n}$. So when $(I',i) \in \mathfrak{J}$, it follows that $I_1' \in \cP$.

We split the sum

\begin{equation}\label{first-split}
    \begin{aligned}
        \sum \sigma(b,l)Z^{I'_b \bI'_l}\nabla_{l}\lambda\nabla_{\bb}\lambda 
    & = \sum\limits_{\{I'\in \mathfrak K\}} \sigma(b,l)Z^{I'_b \bI'_l}\nabla_{l}\lambda\nabla_{\bb}\lambda + \sum\limits_{\{I'\notin \mathfrak K\}} \sigma(b,l)Z^{I'_b \bI'_l}\nabla_{l}\lambda\nabla_{\bb}\lambda\\
    \end{aligned}
\end{equation}

The first term here is controlled using the first term in \eqref{third-order-control2}. For the second term, we use the compatibility equations \eqref{s8} to write

\begin{equation}\label{left-out-third-order} \scalebox{0.95}{$
    \begin{aligned}
      \frac{1}{\lambda^2} \sum\limits_{\{I'\notin \mathfrak K\}} \sigma(b,l)Z^{I'_b \bI'_l}\nabla_{l}\lambda\nabla_{\bb}\lambda &= \sum\limits_{\{I'\notin \mathfrak K\}} \sigma(b,l)Z^{I'_b \bI'_l} (\vp' w_{l} + \psi' u_l)(\vp' w_{\bb} + \psi' u_{\bb})\\
       &\leq (1 + \epsilon_1)(\vp')^2 \sum\limits_{\{I'\notin \mathfrak K\}}\sigma(b,l)Z^{I'_b \bI'_l} w_{l}w_{ \bb} + C(\epsilon_1) (\psi')^2 \\
       &\times\sum\limits_{\{I'\notin \mathfrak K\}}\sigma(b,l)Z^{I'_b \bI'_l} u_{l} u_{\bb} 
    \end{aligned}$}
\end{equation}

\noindent Here for each fixed $I'$, the matrix $(Z^{I'_i \bI'_j})_{i,j}$ is positive-definite and we applied Cauchy's inequality separately for each $I'$.

Then,

\begin{equation}\label{J-split1}
    \begin{aligned}
       \sum\limits_{\{I'\notin \mathfrak K\}} \sigma(b,l)Z^{I'_b \bI'_l}u_{l} u_{\bb}& \leq |\nabla u|^2\sum\limits_{\{I'\notin \mathfrak K\}}\sum\limits_{b\notin I'} Z^{I'_b \bI'_b}
    \end{aligned}
\end{equation}

\

For this term, we consider several cases. If $(I',b)$ is not in $\mathfrak J$, there are several possibilities for $I'_b$, as given below. 

\

{\bf Case 1.1: $1 \in I'$}

\

In this case, $Z^{I'_b \bI'_{b}} = Z^{J'_1 \bJ'_{1}},$ where $J' = (I' \setminus \{ 1 \})\cup \{b\} $, and $1 \notin J'$.

\

{\bf Case 1.2: $I'_b \in ({\cP'})^c $}

\

By Lemma \ref{inversion-lemma}

$$Z^{I'_b \bI'_{b}} \leq {(1+ \epsilon)\delta^{-1}}{\lambda^{-1}} + C(\epsilon){\delta^{-2}}{\lambda^{-2}} \sum Z^{I \bI} .$$

{\bf Case 1.3: $ b\in S$.} 

\

In this case, we have $\lambda_b^2 \geq \tau^2 \lambda^2$.

\

Combining the above three cases, we write

\begin{equation}\label{3casescombine1}
\begin{aligned}
     \sum\limits_{I'\notin \mathfrak K} \sigma(b,l)Z^{I'_b \bI'_l} u_{l} u_{\bb} \leq & (p-1)\sum_{1 \notin J'} Z^{J'_1 \bJ'_1} |\nabla u|^2 + C \left({(1+ \epsilon) \delta^{-1}}{\lambda^{-1}} +  C(\epsilon){\delta^{-2}}{\lambda^{-2}} \sum Z^{I \bI}\right) |\nabla u|^2\\
    & + \tau \lambda \sum_{b \in S} Z^{I'_b \bI'_b}\\
    & \leq 2\tau \lambda\sum_{b \in S} Z^{I'_b \bI'_b} + \frac{C}{A^2}\left( 1 + \sum Z^{I \bI}\right)
\end{aligned}
\end{equation}

Here we used $\lambda \geq C(A,\epsilon,\delta,\tau) (1 +|\nabla u|^2)$, for any known constant $C(A,\epsilon,\delta,\tau)$.

\bigskip

Since $(1+ \epsilon_1) (\vp')^2 \leq \vp''$, using \cref{J-split1,3casescombine1}, equation \eqref{left-out-third-order} becomes

\begin{equation} \label{third-order1}
    \begin{aligned}
          \frac{1}{\lambda^2}  \sum\limits_{\{I'\notin \mathfrak K\}} \sigma(b,l)Z^{I'_b \bI'_l}\nabla_{l}\lambda\nabla_{\bb}\lambda 
       &\leq  \vp'' \sum\limits_{\{I'\notin \mathfrak K\}} \sigma(b,l)Z^{I'_b \bI'_l} w_{l}w_{ \bb}\\
       & + A^2 C(\epsilon_1) \tau \lambda \sum\limits_{i \in S} Z^{I'_i \bI'_i}\;\;+ C_2(1+ \sum Z^{I \bI} )\\
    \end{aligned}
\end{equation}

\

\noindent where $C_2$ is a constant independent of $A$.

From the definition of $\mathcal H$ and \eqref{A-coefficient},

\begin{equation}\label{H-control2}
    \begin{aligned}
        \mathcal{H} \geq& \sigma(i,j) Z^{I'_i \bI'_j} \varphi''   w_{j} w_{ \bi} +  \vp' Z^{I'_i \bI'_i} \lambda_i^2 - \alpha \psi' \sum Z^{I \bI}  - CA\\
    & + \psi'' \sigma(i,j) Z^{I'_i \bI'_j} u_{j}u_{\bi}
    \end{aligned}
\end{equation}

We also know from the $C^{0}$ estimate that $-\psi' \geq C_1 A$, for a positive constant $C_1$ bounded away from zero. Writing $\alpha_1 = \alpha C_1$ and using \eqref{H-control2} with \eqref{third-order1}, we get

\begin{dmath}\label{H-control3}
        \cH -  \frac{1}{\lambda^2} \sum\limits_{\{I' \notin \mathfrak K\}}  \sigma(b,l)Z^{I'_b \bI'_l}\nabla_{l}\lambda\nabla_{\bb}\lambda \geq   (\alpha_1 A - C_2)\sum Z^{I \bI}  - CA+\left(\frac{\vp'}{6} \tau^2\lambda^2- A^2C(\epsilon_1)\tau\lambda\right)\sum_{b \in S} Z^{I_b' \bI_b'}  +\psi'' \sigma(i,j) Z^{I'_i \bI'_j} u_j u_{\bi} +  \vp'' \sum\limits_{I' \in \mathfrak K} \sigma(b,l)Z^{I'_b \bI'_l} w_{l}w_{ \bb} + \frac{\vp'}{2}\sum Z^{I'_i \bI'_i} \lambda_i^2 \geq \frac{\alpha_1 A}{2}  \sum Z^{I \bI}  - CA + \psi'' \sigma(i,j) Z^{I'_i \bI'_j} u_j u_{\bi} +  \vp''  \sum\limits_{I' \in \mathfrak K} \sigma(b,l)Z^{I'_b \bI'_l} w_{l}w_{ \bb} 
       + \frac{\vp'}{2}\sum Z^{I'_i \bI'_i} \lambda_i^2\\
\end{dmath}

\noindent where going to the second inequality we did the following. Recall $\vp' \geq \dfrac{1}{4 \sup(1 + |\nabla u|^2)}$.

\begin{enumerate}
    



    \item $\left(\dfrac{\vp'}{6} \tau^2\lambda^2- A^2C(\epsilon_1)\tau\lambda\right) > 0$ by assuming $\lambda$ larger than $\dfrac{6A^2C(\epsilon_1)}{\tau^2 \vp'}$.

    \ 

    \item $(\alpha_1 A - 2C_2)\geq \dfrac{\alpha_1 A}{2}$, by choosing $A > \dfrac{ (4C_2 + 1)}{\alpha_1}.$

\end{enumerate}

\

Choose $ \epsilon_2$ such that $ \epsilon_2\leq \dfrac{\alpha'}{1 + \alpha'}$ for 
\begin{equation}\label{fixalpha}
    \alpha' = \frac{(2 + \sup |u|)^2}{4 A (1 + \sup |u|)^2}.
\end{equation}

 Now choose $\epsilon_0$, $\tau$ and $\delta$ arbitrarily small so that 
 
\begin{equation}\label{fixconstants}
 (1-\epsilon_0)(1+ \delta_0)^{-1}  \geq ( 1 - \epsilon_2).
\end{equation}

We use the first term of \eqref{third-order-control2} to eliminate $\dfrac{(1 - \epsilon_2) }{ \lambda^2}  \sum\limits_{I' \in \mathfrak K} \sigma(i,j) Z^{I'_i \bI'_j} \nabla_{j}\lambda \nabla_{\bi}\lambda $ from the first term of \eqref{first-split}.

\begin{equation}\label{third-order-control1}
    \begin{aligned}
         \frac{1 }{\lambda} \sum_a Z^{I \bJ} Z^{K \bL} \nabla_{a} Z_{L \bI} \nabla_{\ba} Z_{J \bK} - \dfrac{(1 - \epsilon_2) }{ \lambda^2}  \sum\limits_{I' \in \mathfrak K}  \sigma(i,j) Z^{I'_i \bI'_j} \nabla_{j}\lambda \nabla_{\bi}\lambda & \geq - \frac{C}{\epsilon_0 \lambda^2} \sum Z^{I \bI} 
    \end{aligned}
\end{equation}

Now combining the above two equations $\eqref{third-order-control1}$ and \eqref{H-control3},

\begin{equation}\label{s11.1}
    \begin{aligned}
         \cH + \sum_a\frac{1 }{\lambda} Z^{I \bJ} Z^{K \bL} \nabla_{a} Z_{L \bI} \nabla_{\ba} Z_{J \bK} &- \dfrac{1 }{ \lambda^2}  \sigma(i,j) Z^{I'_i \bI'_j} \nabla_j\lambda \nabla_{\bi}\lambda \geq  -  \epsilon_2 \sum\limits_{I' \in \mathfrak K} \sigma(i,j) Z^{I'_i \bI'_j} \frac{\nabla_j\lambda \nabla_{\bi}\lambda }{\lambda^2}  \\
        &  - CA + \psi'' \sigma(i,j) Z^{I'_i \bI'_j} u_j u_{\bi}  + \varphi''\sum\limits_{I' \in \mathfrak K}\sigma(i,j) Z^{I'_i \bI'_j}    w_{j} w_{ \bi} \\
        &+\frac{\vp'}{2}\sum Z^{I'_i \bI'_i} \lambda_i^2 + \frac{\alpha_1 A}{2}  \sum Z^{I \bI}
    \end{aligned}
\end{equation}

The final step is to control the term $-  \epsilon_2\sum\limits_{I'\in \mathfrak K}\sigma(i,j) Z^{I'_i \bI'_j} \dfrac{\nabla_j\lambda \nabla_{\bi}\lambda }{\lambda^2} $. Then by the compatibility equations \eqref{s8},

\begin{dmath}\label{s31}
 \varphi''\sum\limits_{I' \in \mathfrak K} \sigma(i,j)Z^{I'_i \bI'_j}w_{ j} w_{\bi} \geq  \sum\limits_{I' \in \mathfrak K} \sigma(i,j)Z^{I'_i \bI'_j}\left( \frac{\nabla_j\lambda}{\lambda} + \psi' u_j  \right)\left( \frac{\nabla_{\bi}\lambda}{\lambda} + \psi' u_{\bi}  \right)\geq  \epsilon_2 \sum\limits_{I' \in \mathfrak K} \sigma(i,j)Z^{I'_i \bI'_j}\frac{\nabla_{j}\lambda \nabla_{\bi}\lambda}{\lambda^2} - \frac{\epsilon_2 {\psi'}^2}{1- \epsilon_2}\sum\limits_{I' \in \mathfrak K}\sigma(i,j)Z^{I'_i \bI'_j} u_j u_{\bi}
\end{dmath}

\noindent where in the first line we used $\vp'' \geq {\vp'}^2$. Now $\epsilon_2$ chosen above satisfies

\begin{equation}\label{s32}
   \frac{\epsilon_2 {\psi'}^2}{1- \epsilon_2} < \psi''. 
\end{equation}

Combine equations \eqref{s11.1}, \eqref{s11} with \eqref{s14} to get

\begin{equation}\label{final}
    \begin{aligned}
        0\geq \frac{\vp'}{2}\sum Z^{I'_i \bI'_i} \lambda_i^2 +\frac{\alpha_1 A}{4}  \sum Z^{I \bI}  - CA 
    \end{aligned}
\end{equation}

In the above inequalities we used the fact that 

$$\sum\limits_{I'}\sum_{i,j \notin I'}\sigma(i,j)Z^{I'_i \bI'_j} \xi_j \xi_{\bi} \geq \sum\limits_{I' \in \frak K}\sum_{i,j \notin I'}\sigma(i,j)Z^{I'_i \bI'_j} \xi_j \xi_{\bi},$$

\noindent for any subset $\mathfrak{K} \subset \mathfrak I_{p-1}$ and any vector $\xi$. This holds since for each fixed $I'$, $(Z^{I'_i \bI'_j})_{i ,j}$ is positive-definite. 

\

At the point where $\lambda_1 > \dfrac{\lambda}{2n}\gg 1$, first assume that $\sum\limits_I Z^{I \bI} \leq C $, for some known constant $C$. Since $\det(Z_{I \bJ})$ is bounded, this would imply that all the eigenvalues $\Lambda_I(Z)$ and hence $|Z_{I \bI}|$ are also bounded for all $I$. This is possible only if ${\cP'}^c = \emptyset$. But if ${\cP'}^c = \emptyset$, and if $(I',b) \in \fJ$ then $I'_1 \in {\cP'}^c$ by \eqref{s20}, which is not possible. So we get that $\fJ = \emptyset$ which means that for any $I'_b$ either $b \in S$ and $1 \notin I'_b$, or $1 \in I'_b$. Then

\begin{equation}\label{s35.5}
\begin{aligned}
         \frac{\vp'}{2}\sum Z^{I'_i \bI'_i} \lambda_i^2 & \geq \frac{\vp'}{2} (\sum\limits_{b \in S, 1 \notin I'_b }  Z^{I'_b \bI'_b} \lambda_b^2+ \sum_{\{J':\;1 \notin J'\}} Z^{J'_1 \bJ'_1} \lambda_1^2)\\
         &\geq  \frac{\vp'\tau^2\lambda^2 }{2} \sum\limits_{I} Z^{I \bI}
\end{aligned}
\end{equation}

\noindent where $J' = I'\setminus\{1\} \cup \{b\}$ in the case when $1 \in I'_b$. In this case, $-C A $ is controlled by choosing $\lambda$ large, since $\sum Z^{I \bI}$ is the trace of the inverse which is bounded by the determinant of the inverse, that is prescribed by the RHS of the equation \eqref{main-pde}.

On the contrary, if $\sum\limits_I Z^{I \bI} \gg 1$, then $-C A  $ is controlled by the term $\dfrac{A\alpha}{4} \sum\limits_I Z^{I \bI}$.

This means that the RHS of \eqref{final} is strictly positive, giving a contradiction. So our initial assumption is false and $\lambda \leq C \sup (1 + |\nabla u|^2)$ for some constant $C$ independent of $u$.

\

 Assume that $\cP = \emptyset$ and $\cP'$ is non-empty. Then if $(I',b) \in \fJ$, $I'_1 \in \cP$  as in \eqref{s20} giving a contradiction. Hence the set $\fJ$ must be empty. So by compatibility equation \eqref{s8},

 \begin{equation}
     \begin{aligned}
         \frac{1}{\lambda^2} \sigma(b,l)Z^{I'_b \bI'_l}\nabla_l\lambda\nabla_{\bb}\lambda &=   \sigma(b,l)Z^{I'_b \bI'_l} (\vp' w_{l} + \psi' u_l)(\vp' w_{\bb} + \psi' u_{\bb})\\
       &\leq (1 + \epsilon_1)(\vp')^2 \sigma(b,l)Z^{I'_b \bI'_l} w_{l}w_{ \bb} + C(\epsilon_1) (\psi')^2 \sigma(b,l)Z^{I'_b \bI'_l} u_{l} u_{\bb}
     \end{aligned}
 \end{equation}

 Since $\fJ$ is empty, we always have $(I',b)$, $(I',l) \notin \mathfrak J$. So by \eqref{3casescombine1}

\begin{equation}\label{s36}
\begin{aligned}
    \sigma(b,l)Z^{I'_b \bI'_l} u_{l} u_{\bb}=\sum\limits_{(I',b), (I',l) \notin \mathfrak J} \sigma(b,l)Z^{I'_b \bI'_l} u_{l} u_{\bb} \leq & 2\tau \lambda\sum_{b \in S} Z^{I'_b \bI'_b} + \frac{C}{A^2}(1 + \sum Z^{I \bI})
\end{aligned}
\end{equation}


Hence

\begin{equation}\label{s37}
\begin{aligned}
   \mathcal{H} - &\sigma(i,j)Z^{I'_i \bI'_j}\frac{\nabla_j\lambda \nabla_{\bi}\lambda}{\lambda^2} - C \sum Z^{I \bI}\\
   \geq & \; \vp' \sum  Z^{I'_i \bI'_i} \lambda_i^2  + \psi'  Z^{I'_i \bI'_j}u_{j \bi}- {C}(1+\sum Z^{I \bI})  -A^2C(\epsilon_1) \tau \lambda\sum_{b \in S} Z^{I'_b \bI'_b} \\
> &\;\frac{A\alpha_1}{2}  \sum Z^{I \bI} - CA + \frac{\vp'}{2}  \sum  Z^{I'_i \bI'_i} \lambda_i^2
\end{aligned}
\end{equation}

\

 \noindent which follows by choosing $A > \dfrac{(4C+1)}{\alpha_1}$, and $\lambda>\dfrac{6A^2C(\epsilon_1)}{\tau^2 \vp'}$ as before. Plugging this into \eqref{s14} gives a contradiction, and hence

 $$\lambda \leq C( 1 + \sup |\nabla u|^2),$$

 \noindent for some constant $C$.

 \

 Finally if $\cP' = \emptyset$, then all $I$ satisfies $Z_{I \bI} \geq  \delta \lambda$. Then from the compatibility equations \eqref{s8}

\begin{equation}
    \begin{aligned}
       \sigma(i,j)Z^{I'_i \bI'_j}\frac{\nabla_j\lambda\nabla_{\bi}\lambda}{\lambda^2}& \leq (1+ \epsilon_1) (\vp')^2  \sigma(i,j)Z^{I'_i \bI'_j} w_{j} w_{\bi}
    + C(\epsilon_1) (\psi')^2  (1 + |\nabla u|^2) \sum Z^{I'_i \bI'_i}\\
    & \leq (1+ \epsilon_1) (\vp')^2  \sigma(i,j)Z^{I'_i \bI'_j} w_{j } w_{\bi}+C(\epsilon_1)  (\psi')^2(1 + |\nabla u|^2)\\
    &\times\left(({1+ \epsilon}){\delta^{-1}\lambda^{-1}} + {C}{\delta^{-2}\lambda^{-2}} \sum\limits_I Z^{I \bI} \right),
    \end{aligned}
\end{equation}

\noindent so that

\begin{equation}
\begin{aligned}
   & \mathcal{H} - \sigma(i,j)Z^{I'_i \bI'_j}\frac{\nabla_j\lambda \nabla_{\bi}\lambda}{\lambda^2} - C \sum Z^{I \bI}\\
   &\geq (\varphi''-(1+ \epsilon_1) (\vp')^2)  \sigma(i,j)Z^{I'_i \bI'_j} w_{ j} w_{\bi} + \psi'  Z^{I'_i \bI'_j}u_{j \bi}- {C}\sum Z^{I \bI} + \vp' \sum\limits_{i} Z^{I'_i \bI'_i} \lambda_i^2\\
&  -C(\epsilon_1)  (\psi')^2(1 + |\nabla u|^2)\left(({1+ \epsilon}){\delta^{-1}\lambda^{-1}} + {C}{\delta^{-2}\lambda^{-2}} \sum\limits_I Z^{I \bI} \right) + \psi''\sigma(i,j)Z^{I'_i \bI'_j} u_j u_{\bi} \\
&> \frac{A\alpha_1}{2}  \sum Z^{I \bI} - CA + \frac{\vp'}{2} \sum\limits_{i } Z^{I'_i \bI'_i} \lambda_i^2
\end{aligned}
\end{equation}

 \noindent for $A$ large and $\lambda \gg (\psi')^2(1 + |\nabla u|^2)$. The $-CA$ is controlled as before. This combined with \eqref{s14} gives a contradiction again.


        

\
 
 This completes the proof of second-order estimates.
 









\section{Gradient estimates and existence of solutions}\label{section5}

The blow-up argument for gradient estimates has been used by Chen \cite{Chen2000} for the Monge-Amp\`ere equation and later by Dinew-Kolodziej \cite{DK17} for Hessian equations. It was adapted to the setting of $(n-1)$ plurisubharmonic functions by Tosatti-Weinkove \cite{TW17}. We extend it to the equation under consideration. The important step in this technique is a Liouville-type theorem for some functions defined on $\mathbb C^n$.


This section closely follows the arguments in \cite{TW17}. But we include it to give a self-contained treatment. Assume in this section that $1 < p < n$.

\begin{definition}
    A $C^2$ function $u$ on $\mathbb C^n$ is called $p$-plurisubharmonic if 
      $$\lambda_1 + \lambda_2 + \hdots + \lambda_p > 0$$

      \noindent where $\lambda_1 \leq \lambda_2 \leq \hdots \leq  \lambda_n$ are eigenvalues of the Hessian of $u$.
\end{definition}

In the context of differential forms on $(M, \omega)$, this is equivalent to requiring that $\sqrt{-1} \pbp u \wedge \omega^{p-1}$ is a positive $(p,p)$-form, where if $M = \mathbb C^n$, $\omega$ is the Euclidean metric. In particular, taking trace of this form shows that these functions are subharmonic. 

This notion can be extended to when $u$ not $C^2$. For an upper semi-continuous, locally integrable function $u$ on $\Omega \subset \mathbb C^n$ taking values in $[-\infty, \infty)$, define $p$-PSH as the positivity of the $(p,p)$-form $\pbp u \wedge \omega^{p-1}$ as a current. By Proposition \ref{positive-theorem}, any $u \in C^2$ that is $p$-plurisubharmonic is also $p$-PSH. We will denote $p$-PSH functions on a domain $\Omega$ by $p$-$PSH(\Omega)$.

\begin{definition}
    A $p$-$PSH(\Omega)$ function $u$ is called maximal if for any $\Omega' \Subset \Omega'' \Subset \Omega$ and $v$ any $p$-$PSH(\Omega'')$ function with $v \leq u$ on $\partial \Omega'$,
    
    $$v \leq u$$
    
    \noindent in $\Omega'$.
\end{definition}

Denote by $P(u)$ the current $ \sqrt{-1}\pbp u \wedge \omega^{p-1}$, and for a smooth function $u$, let $P(u)^N = \det{\left(\sqrt{-1}\pbp u \wedge \omega^{p-1}\right)}$, for the Euclidean metric $\omega$ on $\mathbb C^n$.

\begin{lemma}\label{lemma5.1}
    Let $\Omega \subset \mathbb C^n$ be a bounded domain and $u,v: \bar{\Omega} \to \mathbb R$ continuous functions smooth on $\Omega$, with $P(u) \geq 0$ and $P(v) \geq 0$, such that $P(v)^N \geq P(u)^N$ on $\Omega$ and $v \leq u$ on $\partial \Omega$. Then $v \leq u$ on $\Omega$. 
\end{lemma}

\begin{proof}

 By a small perturbation $v_{\epsilon}(z) = v(z) + \epsilon\left( |z|^2 - \sup\limits_{z \in \partial \Omega}{ |z|^2}\right)$, we get that

\begin{equation}
    \begin{aligned}
        P(v_{\epsilon}) > 0, \;\;\;\; P(v_{\epsilon})^N> P(u)^N,
    \end{aligned}
\end{equation}

\noindent in $\Omega$, and $v_{\epsilon} \leq u$ on $\partial \Omega$. Then on $\Omega$,

\begin{equation}
    \begin{aligned}
        0 > P(u)^N - P(v_{\epsilon})^N = \left(\int_{0}^{1} G^{i \bj}(tu + (1-t)v_{\epsilon}) dt\right) (u - v_{\epsilon})_{j \bi} 
    \end{aligned}
\end{equation}

\noindent where $G^{i \bj}$ is the linearization of $P^{N}(u)$, which is positive-definite at $tu + (1-t)v$ for all $t\in [0,1]$ by \eqref{ellipticity}. So by maximum principle, $v_{\epsilon} < u$ on $\Omega$. Let $\epsilon \to 0$ to get $v \leq u$ on $\Omega$ as needed.

\end{proof}

\begin{lemma}\label{lemma5.2}
    Let $u_j$ be continuous, $p$-$PSH$ functions on $\mathbb{C}^n$ that are all maximal. Then if $u_j \to u$ locally uniformly in $\mathbb C^n$, then $u$ is continuous, $p$-$PSH$ and maximal.
\end{lemma}
\begin{proof}

Clearly the limit of positive currents is positive, so $u$ is $p$-PSH. For any $p$-PSH function $v$ defined on some $\Omega \subset \mathbb C^n$, and $\Omega' \Subset \Omega$, assume $v \leq u$ on $\partial \Omega'$. Then $v \leq u_j + \epsilon$ on  $\partial \Omega'$ for some $\epsilon$ and all $j$ large. By maximality of $u_j$, $v \leq u_j+ \epsilon$ on $\Omega'$. Take limit as $\epsilon \to 0$ to get that $v < u$ on $\Omega'$. So $u$ is maximal.

\end{proof}

The following lemma is an extension from the $(n-1)$ case.

\begin{lemma}\label{lemma5.3}
    Let $u$ be a $p$-$PSH$ function on $\mathbb C^n$, $n\geq 3$ that is independent of $z_n$ and maximal. Then $u(z_1, \hdots, z_{n-1},0)$ as a function on $z_1, \hdots, z_n$ is maximal and $(p-1)$-PSH on $\mathbb C^{n-1}$. 
\end{lemma}

\begin{proof}

It is clear that the function is $(p-1)$-PSH on $\mathbb C^{n-1}$, since it is already independent of $z_n$. To prove the maximality, we show that for any bounded domains $\Omega \Subset \Omega' \Subset B_{R}(0)$ for a ball of radius $R$ in $ \mathbb C^{n-1}$, and any continuous $(p-1)$-PSH function $v$ on $\Omega'$, the maximality property is satisfied by $u(z_1,\hdots, z_{n-1}, 0)$. First construct an extension of $v$ to a domain in $\mathbb C^n$,

$$\tilde{\Omega} = \Omega \times \{|z|^n < C_{\epsilon}\}$$

\noindent for a constant $C_{\epsilon}$ to be determined. Define

$$v_{\epsilon} (z_1, \hdots, z_n)= v(z_1, \hdots, z_{n-1}) + \epsilon(|z_1|^2 + \hdots + |z_{n-1}|^2 - |z_n|^2) - R^{2} \epsilon,$$

\noindent which is defined on $\Omega' \times \mathbb C$. Then $P(v_{\epsilon})\geq 0$ since $(p-1)$-PSH functions are in particular $p$-PSH.

 It can be shown that $v_{\epsilon} \leq u$ on $\partial \tilde{\Omega}$ by choosing $C_{\epsilon}$ large enough depending on $\sup |u|$ and $\sup v$. So by maximality of $u$, $v_{\epsilon} \leq u$ on $\tilde{\Omega}$. By restriction to $|z_n| =0$,

 $$ v(z_1, \hdots, z_{n-1}) + \epsilon(|z_1|^2 + \hdots + |z_{n-1}|^2 ) - R^{2} \epsilon \leq u(z_1, \hdots, z_{n-1},0).$$

 Taking limits as $\epsilon \to 0$ now shows that $v \leq u$ on $\Omega$ and hence $u$ is maximal in $\mathbb C^{n-1}$.

\end{proof}

\begin{theorem}[Liouville theorem] \label{liouville-theorem}
A bounded, maximal, $p$-$PSH$ function on $\mathbb C^n$ with bounded gradient is constant.
\end{theorem}

 Let $[\alpha]_r = [\alpha_{I \bJ}]_r dz_I \wedge dz_{\bJ}$ for any differential form, and 

$$[f]_{r}(z) = \frac{1}{|B_r|}\int_{B_r(z)} f \omega^n,$$

\noindent is the average of the function $f$ on the ball $B_r(z)$.

We will show that for a positive function $u$ with $P(u) \geq 0$, if $\left[\left| \dfrac{\p u}{ \p w}\right|^2\right]_r   $ is bounded below by a positive constant for any unit vector $w$, then $[P(u^2)]_r$ is also bounded below.
But this is clear since

\begin{equation}
\begin{aligned}
    P(u^2)  &\geq 2 \sqrt{-1} \partial u \wedge \bpartial u \wedge \omega^{p-1}\\
    &= 2 (p-1)! (\sqrt{-1})^{p^2}\sum_{I'} \sum\limits_{i, j \notin I'} u_i u_{\bj} dz_{I'_i}\wedge dz_{\bI'_j} \\
    &= h^{i \bj}dz_{I'_i} \wedge dz_{\bI'_j}
\end{aligned}
\end{equation}

So in coordinates that diagonalize the Hermitian matrix $[h^{i \bj}]_r$, with eigenvalues $[2\sum\limits_{i \in I}|u_i|^2]_r$,

\begin{equation}\label{change1}
    [P(u^2)]_r \geq 2 (p-1)! (\sqrt{-1})^{p^2} \sum_{I}\sum_{i \in I}[ |u_{i }|^2]_r dz_{I} \wedge dz_{\bI} 
\end{equation}

\noindent is bounded below by a positive form by assumption.

Using this we can complete the proof of Theorem \ref{liouville-theorem}.

\begin{proof}[Proof of Theorem \ref{liouville-theorem}:]

We proceed by induction on the dimension $n$. For $n=2$, the function is plurisubharmonic and the Liouville theorem is known. Let $n\geq 3$, and for contradiction assume that $u$ is a non-constant. Normalize so that $\sup_{\mathbb C^n} u =1$ and $\inf_{\mathbb C^n} u = 0$. Then it can be shown that the following statement and its negation are both false.

There exists a radius $\rho>0$, a sequence of maps on $\mathbb C^n$ of the form $G_k = H_k z + \lambda_k$ with $H_k \in U(n)$ and $\lambda_k \in \mathbb C^n$, and a sequence of radii $r_k \to \infty$ such that

\begin{equation}\label{statement}
    [u^2\circ G_k]_{r_k}(0) + [u\circ G_k]_{\rho}(0) - 2 u\circ G_k(0) \geq \frac{4}{3},
\end{equation}
\noindent and

\begin{equation}
    \lim\limits_{k \to \infty} \int_{B_{r_k}(0)} \left| \frac{\p (u\circ G_k)}{\p z^n}\right|^2 \omega^n =0.
\end{equation}

 Then $u_k = u\circ G_k$ is Lipschitz continuous and maximal. From the assumptions and Lemma \ref{lemma5.2}, there is a subsequence of $u_k$ which converges to some function $v$ which is maximal, $p$-PSH, and Lipschitz with $\sup|\nabla v| \leq C$. It can be shown as in \cite{DK17} that $v$ is independent of $z_n$. Hence by Lemma \ref{lemma5.3} $v$ is $(p-1)$-PSH and maximal in $\mathbb C^{n-1}$. Then by the induction hypothesis, $v$ is a constant. This contradicts \eqref{statement}, since $[u^2\circ G_k]_{r_k}(0) \leq 1$ and $0 \leq v \leq 1$.

Assume that the negation of this statement is true. By Cartan's Lemma for subharmonic functions, for any $z \in \mathbb C^n$

\begin{equation}
    \lim_{r \to \infty} [u]_r(z) = \lim_{r \to \infty} [u^2]_r(z) = 1
\end{equation}

Using this and a change of origin, so that 
$u(0) < \dfrac{1}{12}$, we get that for any $r$ large

\begin{equation}
     [u^2]_{r}(0) + [u]_{\rho}(0) - 2 u(0) > \frac{4}{3}.
\end{equation}

For any $z \in \mathbb C$ and any unit vector $w$, we can choose $G$ such that $G(0) = z$, and $G$ takes $\dfrac{\partial }{\partial z_n}$ to $w$. So by negation of the statement \eqref{statement}, we have that $\left[\left| \dfrac{\p u}{ \p w}\right|^2\right]_r $ is bounded below by a positive constant for some $r$ large and any unit vector $w$.

Consider the set 

$$U = \{2u < [u^2]_r + [u]_{\rho} - \frac{4}{3}\}$$

Then for $\gamma = -c |z|^2$, the set 

$$U_{\gamma} = \{2u < [u^2]_r + [u]_{\rho}+ \gamma - \frac{4}{3}\}$$

\noindent is nonempty, bounded and $U_{\gamma} \Subset U$. Then from \eqref{change1}, for an appropriate choice of $c$, both $2u$ and $[u^2]_r + [u]_{\rho}+ \gamma - \frac{4}{3}$ are continuous, $p$-$PSH$ on $U$, and equal on $\p U_{\gamma}$. By maximality of $2u$, this would imply that

$$2 u \geq [u^2]_r + [u]_{\rho}+ \gamma - \frac{4}{3} $$ 

\noindent inside $U_{\gamma}$, which is a contradiction.

\end{proof}

Using the Liouville theorem, the gradient estimates can be obtained.

\begin{theorem}
    Assume that $u$ is a solution to \eqref{main-pde}. Then there is a constant $C$ that depends only on $|F|_{C^2}$ and other background data such that
    $$\sup\limits_M |\nabla u| \leq C$$
\end{theorem}

\begin{proof}

 The proof is by contradiction. Assume that there is no a priori bound for $|\nabla u|$ in terms of $F$, so that there is a sequence of smooth functions $u_k$ and a bounded set of smooth functions $F_j$ such that 
\begin{equation}\label{blow-up-pde}
    \begin{aligned}
        &\mu_{\Omega_{u_j}} = e^{F_j} \omega^n, \;\;\;\Omega_{u_j} > 0
    \end{aligned}
\end{equation} 

\noindent with 

$$\sup_M |\nabla u_j| = \dfrac{1}{\epsilon_j} \to \infty.$$

Let $x_j \in M$ be the point where $|\nabla u_{j}| = \frac{1}{\epsilon_j}$. Then $\sup_M |u_j| \leq C$, and there is a subsequence for $x_j$ that converges to some $x \in M$. Fix a ball $B_2(0)$ in $\mathbb C^n$ that is identified with a coordinate chart centered at $x$ with $g_{i \bj} = \delta_{i j}$ at $x$. Assume that the $j$ will be large enough so that all $x_j$ are in $B_1(0)$. Define the functions 

$$\hat{u}_j(z) = u_j(\epsilon_j z + x_j),$$

\noindent in $B_{\frac{1}{\epsilon_j}}(0)$. Then the following can be shown.

\begin{equation}
    \begin{aligned}
        &\sup\limits_{B_{\frac{1}{\epsilon_j}}(0)}|\hat u_j|\leq C, \;\;\; \sup\limits_{B_{\frac{1}{\epsilon_j}}(0)} |\nabla \hat u_j| \leq C,\\
        \text{ and }\\
        & |\nabla \hat u_j|(0) = 1, \;\;\; \sup\limits_{B_{\frac{1}{\epsilon_j}}(0)}|\sqrt{-1} \pbp \hat u_j| \leq C'
        \end{aligned}
\end{equation}
\noindent where in the last line we used the estimates for $|\pbp u_j| \leq C(1 + |\nabla u_j|^2)$.

Then by Sobolev embedding theorem, for any $p>1$, $|\hat u|_{W^{2,p}_{\text{loc}}}$, and $|\hat u|_{C^{1, \alpha}_{\text{loc}}}$ are uniformly bounded. So there exists a subsequence of $\hat u_{j}$ that converges in the corresponding topologies, to a function $u$ in $\mathbb C^n$, such that 

\begin{equation}
    \sup_{\mathbb C^n} |u| + \sup_{\mathbb C^n} |\nabla u| \leq C,\;\;\; \nabla u(0) \neq 0
\end{equation}

Now we find the limit of the sequence of equations \eqref{blow-up-pde} under this convergence.

\begin{equation}\label{blow-up1}
    \Omega + \epsilon_j^{-2}\sqrt{-1} \pbp \hat u_j \wedge \omega^{p-1} = \Omega + \sqrt{-1} \pbp u_j \wedge \omega^{p-1} > 0  
\end{equation}

\noindent are positive $(p,p)$-forms, which by Proposition\ref{positive-theorem}, are also positive as currents. Taking weak limit of the expression \eqref{blow-up1}, we get that $P(u) > 0$ as a current. From \eqref{blow-up-pde}

\begin{equation}
 e^{\psi_j} \omega^n = \mu_{\Omega_{u_j}} = \mu_{\Omega_{\epsilon_j^{-2}\hat u_j}}. 
\end{equation}

Taking limit as $ j\to \infty$, we get locally uniform convergence

\begin{equation}\label{current-convergence}
    \det(\epsilon_j^2\Omega + \sqrt{-1}\pbp \hat{u}_j \wedge \omega^{p-1}) \to 0
\end{equation}

\noindent in $\mathbb C^n$.  We show that this implies $u$ is maximal. Then idea is to approximate the $p$-PSH functions by smooth functions, and use Lemma \ref{lemma5.1} to prove a maximality property. This will apply as $P(\epsilon_j^2\Omega + \sqrt{-1}\pbp \hat{u}_j \wedge \omega^{p-1})^N$ can be made arbitrarily small for large $j$.

Let $U \Subset U' \subset \mathbb C^n$, be bounded domains and $v$ a $p$-PSH function in $U$. Assume $v \leq u$ on $\partial U$. Smoothen $v$ by using mollifiers, so that $v_{\delta} = v* \xi_{\delta}$ converges uniformly to $v$ in $U'$ and is $p$-PSH. Choose an $\epsilon>0 $ small, so that for all $j$ large,

$$v_{\delta} \leq \hat u_j + 2 \epsilon' \;\;\text{ on } \partial U.$$

Let $\tilde{v}_{\epsilon} = v_{\delta} + \epsilon |z|^2 - \epsilon C$, for a large constant $C$ so that on $\partial U$, $\tilde v_{\epsilon} \leq \hat u_{j}$. On $U$, $P(\tilde v_{\epsilon}) > 0$, and $P(\tilde v_{\epsilon})^N \geq \epsilon^N$. By \eqref{blow-up1}, we get that on $U'$, $P(\hat u_j + \delta_j |z|^2) \geq 0$, for a sequence $\delta_j \to 0$. We also know from \eqref{current-convergence} that $P(\hat u_j + \delta_j |z|^2)^N \to 0$ as $j \to \infty$. So we have that, for $j$ large enough

\begin{equation}
    P(\tilde v_{\epsilon})^N \geq P(\hat u_j + \delta_j |z|^2)^N
\end{equation}

By Lemma \ref{lemma5.1}, we get that $\tilde v_{\epsilon} \leq \hat u_j + \delta_j |z|^2$ on $U$. Taking limits shows $v \leq u$ on $U$, and hence $u$ is maximal. By Theorem \ref{liouville-theorem}, $u$ is a constant. This contradicts $|\nabla u|(0) \neq 0$.
    
\end{proof}

Using the Evans-Krylov theorem \cite{Evans82} in the setting of complex manifolds, the estimates above can be extended to $C^{2, \alpha}(M)$.  This applies since the equation is uniformly elliptic, and concave in $D^2 u$. Then the standard Schauder estimates for elliptic PDEs improves the regularity to $C^{\infty}$.

We conclude the proof of Theorem \ref{main-theorem} by showing the existence of solutions.
For the continuity method, consider the following family of equations for $t \in [0,1]$.

\begin{equation}\label{continuity-path1}
    \begin{aligned}
         \begin{cases}
            \mu_{\Omega_{u_t}} =  e^{t\psi + b_t } \mu_{\Omega}, \\
            \Omega_{u_t} >0, \;\sup_M u_t  = 0
        \end{cases}
    \end{aligned}
\end{equation}


Here we changed $\omega^n$ in the RHS to $\mu_{\Omega}$ as they only differ by a smooth function.
Fix $k\geq 2$ to be some positive integer, $\alpha>0$ and consider the set 

$$T = \{t \in [0,1]: \; \eqref{continuity-path1} \text{ has a unique $C^{k+2, \alpha}(M)$ solution}\}$$

$T$ is non-empty as $0 \in T$. That $T$ is closed follows from the a priori estimates for $|u_t|_{C^{\infty}(M)}$. To prove the openness of $T$, we apply the inverse function theorem to the linearized operator associated to \eqref{continuity-path1}. Let $t_0$ be a fixed time in $T$, and $u_{t_0}$ be a solution to \eqref{continuity-path1} at $t_0$. Let $G^{i \bj}$ denote the linearization of the operator $\log\mu_{\Omega_v}$ at $u_{t_0}$. By \eqref{ellipticity} $G^{i \bj}$ is a Hermitian metric, and by Gauduchon's theorem \cite{Gauduchon84}, there is a Gauduchon metric $\alpha$ such that 

$$ e^{\sigma } G^{i \bj} \mu_{\Omega_{u_{t_0}}} = \alpha^{i \bj} \alpha^n$$

\noindent for some smooth function $\sigma$, and metric $\alpha$ chosen so that $(\int_M e^{\sigma} \mu_{\Omega_{u_{t_0}}})=1$.

We consider a modified family of equations given by
\begin{equation}\label{continuity-path}
    \begin{aligned}
         \begin{cases}
            \mu_{\Omega_{u_t}} = (\int_M e^{\sigma} \mu_{\Omega_{u_t}}) e^{t\psi + c_t } \mu_{\Omega}, \\
            \Omega_{u_t} >0, \;\sup_M u_t  = 0
        \end{cases}
    \end{aligned}
\end{equation}

\noindent where $c_t$ is a normalization constant defined by

$$c_t = -\log(\int e^{t \psi + \sigma}\mu_{\Omega}) $$ 

Clearly $u_{t_0}$ solves \eqref{continuity-path} at $t=t_0$, with $c_{t_0}=b_{t_0}$. It is enough to show that there is a solution to \eqref{continuity-path} in a neighborhood of $t_0$, since we can set $b_t = \log{(\int_M e^{\sigma} \mu_{\Omega_{u_t}})} +
 c_t$ in \eqref{continuity-path1}.

Define 

$$E_k = \{ v \in C^{k, \alpha}(M):\; \int e^{v + \sigma} \mu_{\Omega_{u_{t_0}}} = 1 \},$$
 \noindent and 

$$P(\Omega) = \{u \in C^{2}(M):\Omega_u \text{ is positive}, \; \text{ and } \int_M u e^{\sigma} \mu_{\Omega} = 0\}.$$

Define the map $MA: C^{k+2,\alpha}(M) \cap P(\Omega) \to E_k$ by 

$$ MA(u) = \log \frac{\mu_{\Omega_u}}{\mu_{\Omega_{u_{t_0}}}} -  \log {\int e^{\sigma}\mu_{\Omega_{u}}}$$

Then the equation \eqref{continuity-path} can be expressed as 

$$MA(u_{t'}) = (t'-t_0) \psi + c_{t'}- c_{t_0}  $$

For any $v$, $MA(v) \in E_k$, and the tangent space at $v$ is given by

$$({TE_k})_{v} = \{\psi\in C^{k, \alpha}(M): \int_M \psi e^{\sigma+ v}\mu_{\Omega_{u_{t_0}}} = 0\}.$$

The tangent space of $C^{k+2,\alpha}(M) \cap P(\Omega)$ at $u_{t_0}$ is a subset of 

$$C^{k+2,\alpha}(M) \cap \{\xi \in C^2(M): \int_M \xi e^{\sigma}\mu_{\Omega} = 0\}$$

\noindent that we denote by $U_{t_0}$.
Consider the operator $\cL: U_{t_0} \to ({TE_k})_{MA(u_{t_0})}$ defined locally by


$$\cL(v) = \sum G^{i \bj} v_{j \bi}$$

\noindent which is the linearization of $MA$ at $u_{t_0}$. Globally this can be written as $\Delta_G v$, the Laplacian with respect to the Hermitian metric $G^{i \bj}$ at $u_{t_0}$. We show that $\cL$ is invertible at $u_{t_0}$. Since this is an elliptic operator, the injectivity follows from maximum principle. For surjectivity, we adapt an idea of Tosatti-Weinkove \cite{TW17}.  Given some $\psi \in ({TE_k})_{MA(u_{t_0})}$, find a unique function $\xi \in U_{t_0}$ such that

$$\Delta_{\alpha} \xi = e^{\sigma}\frac{\mu_{\Omega_{u_{t_0}}}}{\alpha^n}\psi$$

This is possible since the integral of RHS with respect to $\alpha^n$ is zero and the metric $\alpha$ is Gauduchon. Then it follows that 

$$\Delta_G \xi = e^{- \sigma} \frac{\alpha^n}{\mu_{\Omega_{u_{t_0}}}} \Delta_{\alpha} \xi =  \psi,$$

\noindent showing surjectivity. So the mapping $MA$ is a local homeomorphism at $u_{t_0}$. So there is a unique $C^{k+2, \alpha}(M)$ solution in the neighborhood of $t_0$ for \eqref{continuity-path}, and hence for \eqref{continuity-path1}.

This shows that openness of $T$ and hence $1 \in T$, which corresponds to the equation \eqref{main-pde}. The uniqueness of the solution $(u,b)$ follows immediately from the ellipticity and the normalization assumption $\sup_M u = 0$. This ends the proof of Theorem \ref{main-theorem}.

\section{The Geometric flow}\label{section6}

By an abuse of notation we will denote the coefficient function of the volume form $\mu_{\Omega}$ also by $\mu_{\Omega}$. Then $-\sqrt{-1}\pbp \log\mu_{\Omega}$ is independent of the coordinates used and defines a cohomology class. The class 
$$[-\sqrt{-1} \pbp \log \mu_{\Omega} \wedge \omega^{p-1}] \in H^{p,p} (M, \mathbb R)$$

\noindent is also independent of the choice of the positive form $\Omega$. This can be seen from the fact that for any two positive forms $\Omega$ and $\tilde{\Omega}$, the function $\dfrac{\mu_{\Omega}}{\mu_{\tilde{\Omega}}}$ is globally well-defined on $M$. We will denote it by $[\operatorname{Ric_p}(\Omega)]$, since for $\Omega = \omega^{p}$, $-\sqrt{-1} \pbp \log \mu_{\Omega} \wedge \omega^{p-1}$ will be equal to $\operatorname{Ric}(\omega) \wedge \omega^{p-1}$.


 The equation 

\begin{equation}\label{flow}
\begin{aligned}
    &\frac{\p}{ \p t}  \Omega(z, t) =   \frac{pN}{n}\sqrt{-1} \pbp \log{\mu_{\Omega(z, t)}} \wedge \omega^{p-1} ,\\
     &\Omega(z, 0) = \Omega_0(z) \in \Lambda^{p,p}M
    \end{aligned}
\end{equation}

\noindent in terms of cohomology classes says

\begin{equation}
    \frac{d}{d t}[\Omega] = -\frac{pN}{n}[\operatorname{Ric_p}(\tilde{\Omega})],
\end{equation}

\noindent as an equality of sets for some positive form $\tilde{\Omega}$. So similar to K\"ahler-Ricci flow \cite{SW13}, the solution must satisfy

\begin{equation}
    [\Omega] = [\Omega_0] - \frac{pN}{n}t [\operatorname{Ric_p}(\tilde{\Omega})].
\end{equation}

Hence a necessary condition is that the class $[\Omega_0] -\frac{pN}{n} t [\operatorname{Ric_p}(\tilde{\Omega})]$ contains a positive representative. Let 

\begin{equation}\label{existance-time}
\begin{aligned}
    T = \sup\{&t \in (0, \infty): \; \text{$\exists$ a smooth function $\tilde{u}$, and a positive form $\tilde{\Omega}$ such that }\\& \Omega_0 -\frac{pN}{n} t \operatorname{Ric_p}(\tilde{\Omega}) + \sqrt{-1} \pbp \tilde{u} \wedge \omega^{p-1} >0\}.
\end{aligned}
\end{equation}
 
 Choose any $T' < T''< T$, and define

$$\hat{\Omega}_t = \frac{1}{T''}(T'' - t) \Omega_0 + \frac{t}{T''}\left(\Omega_0 - \frac{pN}{n}T'' \operatorname{Ric_p}(\tilde{\Omega}) + \sqrt{-1} \pbp \tilde{u} \wedge \omega^{p-1}\right).$$

\noindent where the function $\tilde{u}$ is obtained from \eqref{existance-time} at $T''$. Clearly $\hat{\Omega}_t \in [\Omega_0] - \frac{pN}{n}t [\operatorname{Ric_p}(\tilde{\Omega})]$. A solution of \eqref{flow} is given by

\begin{equation}\label{flow-form-solution}
    \Omega_u= \hat{\Omega}_t + \sqrt{-1} \pbp u \wedge \omega^{p-1}
\end{equation}

\noindent where $u$ is the solution of the equation

\begin{equation}\label{parabolic-flow}
\begin{aligned}
    &\frac{\partial u}{\partial t} = \frac{pN}{n}\log\frac{\mu_{\Omega_u}}{\theta}, \text{ on } M \times (0, T']\\
    &u|_{t=0}=0
    \end{aligned}
\end{equation}

\noindent where $\theta = e^{\frac{n\tilde{u}}{pNT''}}\mu_{\tilde{\Omega}}$.

\begin{theorem} Any solution $u$ to \eqref{parabolic-flow} in $M \times [0, T']$ satisfies

\begin{enumerate}[(a)]
    \item \label{part1.1} $|{u}|_{C^0} \leq C$.
    \item \label{part2.1} $|\dot{u}|_{C^0} \leq C$.
    \item \label{part3.1} $c\mu_{\tilde{\Omega}} \leq \mu_{\Omega_u} \leq C\mu_{\tilde{\Omega}}.$
\end{enumerate}

\noindent for independent positive constants $c$, $C$.
\end{theorem}

\begin{proof}

Let $Q = u - At$ for a large constant $A>0$ to be determined. Assume $Q$ attains maximum at $(z_0,t_0)$ with $t_0>0$. Choose a coordinate system around $z_0$, such that at $z_0$, $u_{i \bj} = \lambda_i\delta_{i j}$, $g_{i \bj} = \delta_{i j}$, and $dg = 0$. From \eqref{parabolic-flow}

\begin{equation}
    \frac{\partial Q}{\partial t} = \frac{pN}{n}\log\frac{\mu_{\Omega_{Q}}}{\theta} - A \leq \frac{pN}{n}\log\frac{\mu_{\hat{\Omega}_{t_0}}}{\theta} - A < 0
\end{equation}

\noindent for $A$ chosen to be greater than $\frac{pN}{n}\sup\limits_{M \times[0, T']}\log\dfrac{\mu_{\hat{\Omega}_{t}}}{\theta}$. Here we used that 

$$\mu_{\hat{\Omega}_{t_0}} \geq \mu_{\Omega_Q} = \mu_{\Omega_u}$$

\noindent at $z_0$. Indeed at $z_0$, 

$$\sqrt{-1} \pbp Q \wedge \omega^{p-1} = (p-1)! (\sqrt{-1})^{p^2} (\sum\limits_{i \in I} Q_{i \bi}) dz_{I} \wedge dz_{\bI},$$
\noindent  which is negative semi-definite as a matrix. This gives a contradiction so that $t_0 = 0$. For lower bound on $u$, do the similar argument with $Q = u + At$. 

For $\dot{u}$, we take $Q = \dot{u} - Au$. With similar assumptions as above, at the point of maximum $z_0$ of $Q$, 

\begin{equation}
    \begin{aligned}
        0 &\leq \partial_t Q - \sigma(i,j)Z^{I'_i \bI'_j} Q_{j \bi} \\
        & = Z^{I \bJ} \frac{\partial }{\partial t} (\hat{\Omega}_t)_{J \bI}- A \dot{u} + A Z^{I'_i \bI'_i} u_{i \bi} < AN- A \dot{u} 
    \end{aligned}
\end{equation}
\noindent where the last inequality comes from choosing $A \alpha > \sup|-\frac{pN}{n} \operatorname{Ric_p}(\tilde{\Omega})+ \frac{1}{T''}\sqrt{-1} \pbp \tilde u \wedge \omega^{p-1}| + 1$, and \eqref{A-coefficient}. This gives $\sup \dot{u} \leq C$ for $C$ depending on $\sup|u|$. For lower bound, similar to \cite{SW13,ST16} take

$$Q = (T'' - t) \dot{u} + u + Nt,$$

\noindent so that at the point of minimum of $Q$,

\begin{equation}
\begin{aligned}
    0 &\geq  \partial_t Q - \sigma(i,j)Z^{I'_i \bI'_j} Q_{j \bi} \\
    &= (T'' - t) Z^{I \bJ} \frac{\partial (\hat{\Omega}_t)_{J \bI}}{\partial t} + N - Z^{I'_i \bI'_i} u_{i \bi}\\
    &= (T'' - t) Z^{I \bJ} \frac{\partial (\hat{\Omega}_t)_{J \bI}}{\partial t} + Z^{I \bJ} (\hat{\Omega}_t)_{J \bI}\\
    &= Z^{I \bJ}( \Omega_0 - \frac{pN}{n}T'' \operatorname{Ric_p}(\tilde{\Omega}) + \sqrt{-1} \pbp \tilde u \wedge \omega^{p-1} )_{J \bI}> 0
\end{aligned}
\end{equation}

\noindent by \eqref{trace-Z} and \eqref{existance-time}. This gives the required lower bound for $\dot u$ in $[0, T']$. Part \ref{part3.1} follows from $|\dot u |\leq C$ and \eqref{parabolic-flow}.

\end{proof}





To prove Theorem \ref{existence-time-theorem}, we need to show that the flow exists in $(0, T]$. This involves finding a priori $C^{2, \alpha}$ estimate for $u$ in $ (0, T']$ for $T' < T$. We can derive $C^2$ estimates for the solution that depends quadratically on the gradient of $u$ by following the same argument as in the elliptic case. The argument has been outlined for the normalized flow in Theorem \ref{normalized-estimate-theorem}. A parabolic version of the blow-up argument from Section \ref{section5} then shows that $|u|_{C^{1}} \leq C$. Finally using standard parabolic theory, we obtain that $|u|_{C^{\infty}}$ is uniformly bounded.

To extend the flow assume for contradiction that the flow exists only till $T'<T$. Then by short-time existence result of parabolic equations, one can extend the flow past $T'$ by using the initial smooth data at $T'$ obtained from the estimates above. This would give a contradiction so that the flow exists till $T$.

The $T$ above need not be maximal for the evolution of forms, since a solution to \eqref{flow} need not necessarily be of the form of \eqref{flow-form-solution}. But if there is a positive form $\tilde{\Omega}$, such that $\operatorname{Ric_p}(\tilde{\Omega}) \leq 0$, and $\Omega_0$ is positive, then the solution exists for all time. In this setting, we would like to discuss the convergence. If $\operatorname{Ric_p}(\tilde{\Omega})=0$, then the convergence of the normalized solution

$$\tilde{u} = \frac{1}{V} \int_{M} u \omega^n$$

\noindent can be proven by a Harnack inequaltiy for $u_t$ and exponential decay of the oscillation of the normalized solution similar to \cite{George21, Cao85, Gill11}. But in the case $\operatorname{Ric_p}(\tilde{\Omega})> 0 $, we need to reparametrize the flow. This will be done for a modified version of the flow given by \eqref{normalized-form-flow}.



For a positive form $\tilde{\Omega}$, we consider 

\begin{equation}\label{flow3}
\begin{aligned}
   &\frac{\partial \Omega}{\partial t} = \frac{Np}{n}\sqrt{-1} \pbp \log{\frac{ \mu_{\Omega}}{\mu_{\tilde{\Omega}}}}\wedge \omega^{p-1} + \Omega_0\\
   &\Omega|_{t=0} = \Omega_0 \in H^{p,p}(M, \mathbb R)
\end{aligned}
\end{equation}

A solution can be obtained as

\begin{equation}
    \Omega_u = (1 + t) \Omega_0 + \sqrt{-1} \pbp u \wedge \omega^{p-1}
\end{equation}

\noindent for a solution $u$ of 

\begin{equation}\label{parabolic-ma-eqn1}
    \frac{\partial u}{\partial t} = 
    \frac{Np}{n}\log\frac{\mu_{\Omega_u}}{\mu_{\tilde{\Omega}}} 
\end{equation}

\noindent with $u(0) = 0$. $\Omega_0$ being positive will imply the long-time existence of the flow. Similar to K\"ahler-Ricci flow, we use the reparametrization

\begin{equation}
\hat{\Omega}(s) = (s+1)\Omega(\log(s+1))
\end{equation}

\noindent so that if $\Omega$ solves

\begin{equation}\label{normalized-form-flow}
\begin{aligned}
 &\frac{\partial\Omega}{\partial t} =  \frac{Np}{n}\sqrt{-1} \pbp \log{\frac{ \mu_{\Omega}} {\mu_{\tilde{\Omega}}}}\wedge \omega^{p-1}+ \Omega_0 - \Omega\\
   &\Omega|_{t=0} = \Omega_0 \in H^{p,p}(M, \mathbb R)
\end{aligned}
\end{equation}

\noindent $\hat{\Omega}$ is a solution to \eqref{flow3}. A solution to \eqref{normalized-form-flow} is given by 

\begin{equation}
    \Omega_u = \Omega_0 + \sqrt{-1} \pbp u \wedge \omega^{p-1}
\end{equation}

\noindent where $u$ is the solution to the normalized Monge-Amp\`ere equation

\begin{equation}\label{normalized-ma-eqn}
\begin{aligned}
    & \frac{\partial u}{\partial t}  = 
    \frac{pN}{n}\log\frac{\mu_{\Omega_u}}{\mu_{\tilde{\Omega}}} - u, \text{ on } M\\
   & u(0) = 0.
\end{aligned}
\end{equation}

Now we show the long-time existence of \eqref{normalized-ma-eqn} and smooth convergence to a solution $u_{\infty}$, assuming that $\Omega_0>0$. In particular, this would imply that $[\Omega_{\infty}] = [\Omega_0]$.

\begin{theorem}\label{normalized-estimate-theorem}
   Let $u$ be a smooth solution to \eqref{normalized-ma-eqn}. Then

    \begin{enumerate}[(a)]
        \item \label{part1} $ |u|_{C^0}\leq C$.
        \item \label{part2}$|\dot u|_{C^0} \leq C$.
        \item \label{part3} $c\mu_{\tilde{\Omega}} \leq \mu_{\Omega_u} \leq C\mu_{\tilde{\Omega}}.$

        \ 
        
        \item \label{part4} $\sup\limits_{M \times [0, t_1]}|\sqrt{-1} \pbp u|_{\omega} \leq C (1 + \sup\limits_{M \times [0, t_1]}|\nabla u|^2) $.
    \end{enumerate}
    \noindent for any $t_1$ in the existence interval and positive constants $c$ and $C$ independent of $t_1$.
\end{theorem}

\begin{proof}
  We compute immediately that
    \begin{equation} \label{normalized-proof-1}
        \frac{\partial}{\partial t} (e^t \dot{u}) = \sigma(i,j)Z^{I'_i \bI'_j} (e^t \dot{u})_{j \bi}
    \end{equation}
By maximum principle, it follows that

\begin{equation}\label{normalized-proof-2}
    |\dot{u}| \leq C e^{-t}
\end{equation}

To be precise, the above is a consequence of classical maximum principle applied to $v_{\epsilon} = e^t\dot u - \epsilon t$. At the point of maximum of $v_{\epsilon}$ in $(0, T]$, $(\partial_t - \sigma(i,j)Z^{I'_i \bI'_j}\nabla_j \nabla_{\bi}) v_{\epsilon} <0$ by \eqref{normalized-proof-1}, which is a contradiction. It follows by integration of \eqref{normalized-proof-2} that 

\begin{equation}
    |u(z, s) -  u(z, t)| \leq C(e^{-t} - e^{-s})
\end{equation}

Hence $u$ converges uniformly to a continuous function $u_{\infty}$ as $t \to \infty$. Part \ref{part1}, \ref{part2} and \ref{part3} follows from this combined with \eqref{normalized-ma-eqn}. For part \ref{part4}, we need to derive $C^2$ estimates for the normalized parabolic equation, similar to the elliptic case. But since the argument is essentially the same, we only outline the idea.

Consider the same test function 

$$Q = \log(\lambda) + \varphi(|\nabla u|^2) + \psi(u).$$

Assume as before $\lambda \gg \sup (1 + |\nabla u|^2)$. Then by applying $\dfrac{\partial}{\partial t} -\sigma(i,j)Z^{I'_i \bI'_j} \nabla_j \nabla_{\bi}$ to $Q$, we get at the point $z_0$ where $Q$ attains maximum and where the Hessian is diagonalized with $g_{i \bj} = \delta_{ij}$, $dg_{i \bj} =0$,

\begin{equation}\label{norm-c21}
\begin{aligned}
0 &\geq \sum_a\sigma(i,j)Z^{I'_i \bI'_j} \frac{u_{a \ba j \bi}}{\lambda}-\sigma(i,j)Z^{I'_i \bI'_j}\frac{\nabla_{j}\lambda\nabla_{\bi}\lambda}{\lambda^2}+  2 \vp ' \sigma(i,j)Z^{I'_i \bI'_j} \sum_k\fRe{ (u_{\bk} u_{k j \bi})}  + \mathcal{H}\\
&- \frac{(\lambda)_t}{\lambda} - \vp' (u_k u_{\bk t} + u_{\bk} u_{k t}) - \psi' u_t - C \phi' |\nabla u|^2 \sum Z^{I \bI}
\end{aligned}
\end{equation}

Differentiating \eqref{normalized-ma-eqn} by $\Delta$ gives

\begin{equation}\label{norm-c22}
   \lambda_t = - \lambda + \frac{pN}{n}\Delta \log{\frac{\mu_{\Omega_u}}{\mu_{\tilde{\Omega}}}} + (n-p+1)C_0
\end{equation}

Combining \eqref{norm-c21} and \eqref{norm-c22},

\begin{equation}\label{norm-c23}
\begin{aligned}
    0 &\geq \sum_a\sigma(i,j)Z^{I'_i \bI'_j} \frac{u_{a \ba j \bi}}{\lambda}-\sigma(i,j)Z^{I'_i \bI'_j}\frac{\nabla_{j}\lambda\nabla_{\bi}\lambda}{\lambda^2} +  2 \vp ' \sigma(i,j)Z^{I'_i \bI'_j} \fRe{ (u_{\bk} u_{k j \bi})}  + \mathcal{H}\\
& - 2 \vp'  \fRe( u_{\bk} u_{k t}) - \psi' u_t   +1 - \frac{(n-p+1)C_0}{\lambda}  - \frac{pN}{n\lambda}\Delta \log{\frac{\mu_{\Omega_u}}{\mu_{\tilde{\Omega}}}} - C(1+ \sum Z^{I \bI})
\end{aligned}
\end{equation}

Also

\begin{equation}\label{norm-c24}
\begin{aligned}
    |u_{kt } u_{\bk}| &= |\sigma(i,j)Z^{I'_i \bI'_j} u_{j \bi k} u_{\bk} + Z^{I \bJ} u_{\bk} \nabla_k X_{J \bI} - |u_k|^2| \\
    &\geq - C(1 + |\nabla u|^2) \sum Z^{I \bI} + 2 \fRe(\sigma(i,j)Z^{I'_i \bI'_j} u_{k j \bi } u_{\bk} )
\end{aligned}
\end{equation}

Combine \eqref{norm-c23} and \eqref{norm-c24} and use $\vp'  \geq \dfrac{1}{ 8 \sup ( 1 + |\nabla u|^2)}$ to get

\begin{equation}\label{parabolic-C^2}
    \begin{aligned}
         0 &\geq \sum_a\sigma(i,j)Z^{I'_i \bI'_j} \frac{u_{a \ba j \bi}}{\lambda}-\sigma(i,j)Z^{I'_i \bI'_j}\frac{\nabla_{j}\lambda\nabla_{\bi}\lambda}{\lambda^2}   + \mathcal{H} - C \sum Z^{I \bI}\\
& \;\;- \frac{pN}{n\lambda}\Delta \log{\frac{\mu_{\Omega_u}}{\mu_{\tilde{\Omega}}}} - \psi' u_t
    \end{aligned}
\end{equation}

The term $|- \psi' u_t| \leq C A$ can be controlled as before, by the argument involving \eqref{s35.5}. The rest of the proof proceeds in the same way as the elliptic case to get a contradiction that the RHS in \eqref{parabolic-C^2} is positive. This gives the required estimate.
    
\end{proof}

For gradient estimates, the blow-up argument above can be adapted to the parabolic setting. We omit the argument here. This gives that $|\sqrt{-1} \pbp u|_{C^2} \leq C$, and by standard parabolic theory, this can be extended to $|u|_{C^{\infty}} \leq C$.

As the proof of Theorem \ref{normalized-estimate-theorem} shows there is a continuous function $u_{\infty}$ to which $u$ converges to. This convergence is smooth as can be seen from the following argument. Assume for contradiction that it is not smooth. So there exits an integer $k$, an $\epsilon>0$, and a sequence ${t_l}$ such that

\begin{equation}\label{smooth-convergence}
    |u_{t_l} - u_{\infty}|_{C^k} > \epsilon,
\end{equation}

\noindent for all $l$ large. Using the Arzela-Ascoli theorem and the $C^{\infty}$ estimates, we can find a subsequence $t_{l_i}$ such that it converges in $C^k$ to a $u'$. But by \eqref{smooth-convergence}, $|u' - u_{\infty}|_{C^k} > \epsilon$. This means $u' \neq u_{\infty}$, which contradicts the fact that both are also equal to the unique $C^0$ limit of $u$.

Now taking the limit $t \to \infty$ of the flow \eqref{normalized-form-flow}, we get

\begin{equation}
   \Omega_{\infty} = \Omega_0 + \frac{pN}{n}\sqrt{-1} \pbp \log{\frac{ \mu_{\Omega_{\infty}}}{\mu_{\tilde{\Omega}}}}\wedge \omega^{p-1}.
\end{equation}


By taking limit of \eqref{normalized-ma-eqn}, we also get

\begin{equation*}
    \operatorname{Ric}(\Omega_{\infty}) = \operatorname{Ric}(\tilde{\Omega})  -\frac{n}{pN}\sqrt{-1} \pbp u_{\infty}.
\end{equation*}

\vspace{1cm}

{\bf Data Availability Statement:} DAS is not applicable to this article as no datasets were generated or analyzed during the current study.

\

{\bf Conflict of Interest:}
There is no conflict of interest.

\bigskip
\small

\bibliographystyle{plain}
\bibliography{references}

\end{document}